# ASYMPTOTIC INFERENCE IN SOME HETEROSCEDASTIC REGRESSION MODELS WITH LONG MEMORY DESIGN AND ERRORS

By Hongwen Guo and Hira L. Koul

*Michigan State University*

This paper discusses asymptotic distributions of various estimators of the underlying parameters in some regression models with long memory (LM) Gaussian design and nonparametric heteroscedastic LM moving average errors. In the simple linear regression model, the first-order asymptotic distribution of the least square estimator of the slope parameter is observed to be degenerate. However, in the second order, this estimator is $n^{1/2}$-consistent and asymptotically normal for $h + H < 3/2$; nonnormal otherwise, where $h$ and $H$ are LM parameters of design and error processes, respectively. The finite-dimensional asymptotic distributions of a class of kernel type estimators of the conditional variance function $\sigma^2(x)$ in a more general heteroscedastic regression model are found to be normal whenever $H < (1 + h)/2$, and non-normal otherwise. In addition, in this general model, $\log(n)$-consistency of the local Whittle estimator of $H$ based on pseudo residuals and consistency of a cross validation type estimator of $\sigma^2(x)$ are established. All of these findings are then used to propose a lack-of-fit test of a parametric regression model, with an application to some currency exchange rate data which exhibit LM.

**1. Introduction.** This paper discusses asymptotic distributions of some estimators of the underlying parameters in some heteroscedastic regression models with LM in design and errors. This is of interest partly for the purpose of regression model diagnostics and partly for the sake of some large sample inference in these models. Regression models with LM in both design and error variables are useful when LM in the given response variable is not fully explained by LM in the given design variable; see [26].

For the sake of clarity of exposition, we first focus on a simple linear regression model where one observes a strictly stationary bivariate process $(X_t, Y_t)$, $t \in \mathbb{Z} := \{0, \pm 1, \ldots\}$, both having finite and positive variances and obeying the model

(1.1) $\quad Y_t = \beta_0 + \beta_1 X_t + \sigma(X_t) u_t, \qquad \text{for some } (\beta_0, \beta_1) \in \mathbb{R}^2,$

(1.2)
$$u_t := \sum_{j=0}^{\infty} b_j \varepsilon_{t-j}, \qquad b_j \sim C j^{-(3/2-H)},$$
$$\text{as } j \to \infty, \text{for some } \tfrac{1}{2} < H < 1.$$

Here, $\varepsilon_t$ are standardized i.i.d. r.v.'s, independent of the $X_t$-process and the constant $C$ is such that $\sum_{j=0}^{\infty} b_j^2 = 1$. Under this set up, $\sigma^2(x) = \text{Var}(Y_0 | X_0 = x)$, $x \in \mathbb{R}$, and $E\sigma^2(X_0) < \infty$.

For a stationary second-order process $\xi_t$, $t \in \mathbb{Z}$, let $f_\xi$ ($\gamma_\xi$) denote its spectral density (auto-covariance function). We also assume that $\{X_t\}$ is a Gaussian process with mean $\mu$, variance $\gamma^2 := \gamma_X(0)$, and

(1.3) $\quad \gamma_X(k) \sim \dfrac{G_X \theta(h)}{k^{2(1-h)}} \qquad \text{as } k \to \infty, \text{ for some } 1/2 < h < 1,$

where $\theta(h) := 2\Gamma(2 - 2h) \cos(\pi(1 - h))$ and $G_X > 0$ is a constant. The sequence $b_j$ is also assumed to satisfy $b_j \to 0$, as $j \to \infty$ and for some $a < \infty$, $b_{j+1} \leq b_j(1 + j^{-1}a)$, for all sufficiently large $j$. This condition, for example, is satisfied by FARIMA$(0, H - 1/2, 0)$ model where $b_j = \Gamma(j + H - 1/2)/\Gamma(j + 1)\Gamma(H - 1/2)$. As pointed out in [24], page 1632, under this condition,

(1.4)
$$f_u(\lambda) \sim G_u \lambda^{1-2H}, \qquad \lambda \to 0_+;$$
$$\gamma_u(k) \sim G_u \theta(H) k^{-2(1-H)}, \qquad k \to \infty,$$

where $G_u$ is a positive constant.

Several authors have discussed regression models with LM errors when $\sigma(x) \equiv c$, a constant. The asymptotic distributions of the least squares estimator (LSE) and M- and R-estimators in nonrandom design linear regression models with LM errors are established in [12, 17, 18, 33, 34] and for nonlinear regression models in [19]. The asymptotic distribution of the generalized LSE (GLSE) in certain polynomial regression models is discussed in [6] and [16]. The errors in [6] are assumed to be LM Gaussian, while in [16] a function of a long memory moving average (LMMA) process.

In the context of homoscedastic multiple linear regression models with LM in both covariates and errors and when the error process has a known parametric spectral density, the GLSE of the slope parameter vector is known to be $n^{1/2}$-consistent and asymptotically normal with the Gauss–Markov variance; see [26]. This result is adapted in [15] to the case where the error spectral density is semi-parametric as in (1.4). A crucial result needed here is the availability of a preliminary $n^{1/2}$-consistent estimator of the slope parameter vector. In [26] it was also noted that the LSE is $n^{1/2}$-consistent for certain values of $H, h$. In a simulation study reported in [15] it was found that the adaptive estimator where LSE was used as a preliminary estimator had smaller MSE compared to the one where GLSE was used. This partly motivates the need to understand asymptotic behavior of the LSE in the current set up, for were one to carry out an analogous adaptation program here, even to get started one would need a preliminary $n^{1/2}$-consistent estimator of $\beta_1$ in order to estimate $\sigma(x)$ and $f_u(\lambda)$ nonparametrically. And at least for those values of $H$ and $h$ for which the LSE is $n^{1/2}$-consistent, its use in constructing such adaptive estimators would be justified. Other reasons are to be able to estimate $G_u$ and $H$ and to understand asymptotic behavior of some lack-of-fit tests for fitting a heteroscedastic linear regression model. Currently there is a void on this topic in the literature which this paper is attempting to fill. Because of its simplicity, it is desirable to use the LSE for these purposes.

Section 2 discusses asymptotic distribution of the LSE $\hat{\beta} := (\hat{\beta}_0, \hat{\beta}_1)'$ of $\beta := (\beta_0, \beta_1)'$ in the models (1.1) and (1.2). The weak limit of $n^{1-H}(\hat{\beta} - \beta)$ is shown to be a bivariate normal distribution, for all $1/2 < H, h < 1$. But if $\mu = EX = 0$ and $\sigma(x)$ is an even function, then this asymptotic distribution of $\hat{\beta}_1$ is degenerate. In this case we further obtain that if $H + h < 3/2$, $n^{1/2}(\hat{\beta}_1 - \beta_1)$ converges weakly to a normal r.v. On the other hand, in the case $H \wedge h := \min(H, h) > 3/4$, and even when both $u_t$ and $X_t$ are Gaussian, $\hat{\beta}_1$ has a nonnormal weak limit with the normalization $n^{2-H-h}$.

To implement the proposed lack-of-fit test for fitting a regression model or to carry out some inference about $\beta$ and $\sigma^2(x)$, one needs consistent estimators of $\sigma(x)$, $G_u$ and a $\ln(n)$-consistent estimator of $H$. Section 3 derives asymptotic distributions of a class of kernel type estimators $\hat{\sigma}^2(x)$ of $\sigma^2(x)$ in the regression model

(1.5) $$Y_t = \beta' r(X_t) + \sigma(X_t) u_t, \qquad \beta \in \mathbb{R}^q,$$

where $r(x)$ is a vector of some known $q$ functions and the rest of the entities are as in (1.2) and (1.3). It is proved that when $H < (1+h)/2$, the finite-dimensional distributions of $n^{1-h}(\hat{\sigma}^2 - \sigma^2)$ converge weakly to $k$-variate normal distributions, while for $H > (1+h)/2$, the weak limit of $n^{2-2H}(\hat{\sigma}^2(x) - \sigma^2(x))$ is nonnormal.

Using the approach in [25], the local Whittle estimator of $H$ based on the pseudo residuals $Y_t - \tilde{\beta}' r(X_t)$ in the model (1.5) is shown to be $\log(n)$-consistent, where $\tilde{\beta}$ is the LSE. This is unlike the case of nonparametric heteroscedastic regression model with $X_t = t/n$, $1 \leq t \leq n$, and LMMA errors, where it is necessary to base estimators of $H$ on the standardized residuals; see [11].

An important inference problem is to assess the accuracy of an assumed regression model. Let $(X, Y)$ denote a copy of $(X_0, Y_0)$ and $\mu(x) := E(Y|X =$ $x$). Consider the problem of testing $H_0: \mu(x) = \beta' r(x)$, for some $\beta \in \mathbb{R}^q$ and for all $x \in \mathbb{R}$, against the alternative that $H_0$ is not true. In the 1990's, several authors found that tests of $H_0$ based on the marked empirical process

$$\widetilde{\mathcal{V}}_n(x) = \sum_{t=1}^n (Y_t - \tilde{\beta}' r(X_t)) I(X_t \leq x), \qquad x \in \bar{\mathbb{R}} := [-\infty, \infty],$$

have desirable level and power properties against a broad class of alternatives; see, for example, [[1], [21], [29], [30], [31]], among others. See [1], pages 132–134 of [13], and [29] for more motivation about using this process for lack-of-fit testing. In the presence of long memory in design and/or errors and when $\sigma(x) \equiv c$, some tests based on this process have been studied in [20].

Under the current set up, Theorem 5.1 below proves that, under $H_0$, $n^{-H} \widetilde{\mathcal{V}}_n(x)$ converges weakly to $J_\sigma(x) \psi_1 Z$, in $\mathcal{D}(\bar{\mathbb{R}})$ and uniform metric, where $Z$ is a $N(0,1)$ r.v., $\psi_1^2 := G_u \theta(H)/H(2H-1)$, $J_\sigma(x) := E[\sigma(X) - E(\sigma(X) r(X))' A^{-1} r(X)] \times I(X \leq x)$, and $A^{-1}$ is assumed to exist. Thus, to use $\widetilde{\mathcal{V}}_n$ for testing $H_0$, we need a uniformly consistent estimator of $J_\sigma(x)$ and a consistent estimator of $\psi_1$. A uniformly consistent estimator of $J_\sigma$, under $H_0$, based on the leave-one-observation-out estimator of $\sigma(x)$, is given in Section 5. The regular kernel type estimator is not useful here because of the unstable behavior of $\hat{\sigma}^2(X_t)$. The estimators of $G_u$ and $H$ constructed in Section 4 are used to provide a consistent estimator of $\psi_1$ under $H_0$.

Section 6 includes a finite sample simulation and an application to some monthly currency exchange rate data that exhibits long memory. The last section is the Appendix consisting of some proofs and necessary lemmas.

In this paper all limits are taken as $n \to \infty$, unless specified otherwise. For any two sequences of real numbers, $a_n \sim b_n$ means that $a_n/b_n \to 1$, $\to_d$ stands for the convergence in distribution of a sequence of r.v.'s, while $\Longrightarrow$ denotes the weak convergence of a sequence of stochastic processes, and $u_p(1)$ denotes a sequence of stochastic processes that tends to zero uniformly over its time domain, in probability. Henceforth, the independence of $X_t$ and $u_t$ processes is used without mention.

**2. Asymptotics of the LSE.** This section discusses asymptotic distribution of the LSE in the model (1.1)–(1.3). For this purpose, we need the following result. Let $\nu$ be a real valued function on $\mathbb{R}$ with $E\nu^2(X) < \infty$. Set $\nu_0 := E\nu(X)$. By (A.20) below, there is a $C < \infty$ such that $|E(\nu(X_0) - \nu_0) u_0 (\nu(X_t) - \nu_0) u_t| \leq C t^{-4+2H+2h}$, for all sufficiently large and positive $t$. Hence, $\forall 1/2 < h, H < 1$,

$$n^{-H} \sum_{i=1}^n \nu(X_i) u_i = \nu_0 n^{-H} \sum_{i=1}^n u_i + n^{-H} \sum_{i=1}^n (\nu(X_i) - \nu_0) u_i$$

(2.1)
$$= \nu_0 n^{-H} \sum_{i=1}^n u_i + o_p(1).$$

Next, let $D(a) := \theta(a)/a(2a-1)$, $1/2 < a < 1$, and $Z_1, Z_2$ be two independent r.v.'s, $Z_j$ having $N(0, \psi_j)$ distribution, $j = 1, 2$, where $\psi_1^2 = G_u D(H)$, $\psi_2^2 = G_X D(h)$. Let $\sigma_0 := E\sigma(X)$, $J := EX\sigma(X)$, $\gamma^2 \Gamma_1 := (J - \mu\sigma_0) = \text{Cov}(X, \sigma(X))$ and $\Gamma_0 := \sigma_0 - [(J - \mu\sigma_0)\mu/\gamma^2]$. From [8], we obtain

(2.2) $$Z_{n1} := n^{-H} \sum_{i=1}^n u_i \xrightarrow{d} Z_1, \qquad n^{-h} \sum_{i=1}^n (X_i - \mu)/\gamma \xrightarrow{d} Z_2.$$

Now, let $\sigma_t := \sigma(X_t)$, $e_t := \sigma_t u_t$, $n\bar{X} := \sum_{t=1}^n X_t$, $n\bar{u} := \sum_{t=1}^n u_t$, $n\bar{Y} := \sum_{t=1}^n Y_t$, $n\bar{e} := \sum_{t=1}^n e_t$, $ns^2 := \sum_{t=1}^n (X_t - \bar{X})^2$. Recall that the LSE satisfies

(2.3) $$\hat{\beta}_1 - \beta_1 = \frac{1}{s^2}\left(n^{-1}\sum_{t=1}^n X_t e_t - \bar{X}\bar{e}\right), \qquad \hat{\beta}_0 - \beta_0 := \bar{e} - \bar{X}(\hat{\beta}_1 - \beta_1).$$

From (2.1) and (2.2) applied to $\nu(x) \equiv \sigma(x)$ we obtain that $\bar{e} = O_p(n^{H-1})$, $\bar{X} - \mu = O_p(n^{h-1})$. These facts and (2.3) yield that $n^{1-H}(\hat{\beta}_0 - \beta_0) = \Gamma_0 Z_{n1} + o_p(1)$, $n^{1-H}(\hat{\beta}_1 - \beta_1) = \Gamma_1 Z_{n1} + o_p(1)$ and hence, the following:

LEMMA 2.1. *Under* (1.1)–(1.3), $n^{1-H}(\hat{\beta}_0 - \beta_0, \hat{\beta}_1 - \beta_1) \rightarrow_d (\Gamma_0, \Gamma_1)Z_1$.

But $\Gamma_1 = 0$ if either $\sigma(x) \equiv c$, a constant or $\mu = 0$ and $\sigma(x)$ is an even function of $x$. Since in these cases the weak limit of $n^{1-H}(\hat{\beta}_1 - \beta_1)$ is degenerate at zero, it is pertinent to investigate higher-order approximation to the distribution of $\hat{\beta}_1$. The former case has been discussed in [20]. We shall next discuss the second-order result under the following:

ASSUMPTION 1. $\sigma(x)$ is an even function of $x \in \mathbb{R}$ and $\mu := E(X) = 0$.

Let $a(z) := \int_0^1 u^{z-3/2}(1-u)^{1-2z}\, du$, $1/2 < z < 1$, $\tilde{C} := (\frac{G_u G_X}{a(H)a(h)})^{1/2}$ and

(2.4) $$\mathcal{Z}_2 := \tilde{C} \int \int_0^1 [(s-x_1)^{-(3-2H)/2}(s-x_2)^{-(3-2h)/2}]$$
$$\times I(x_1 < s, x_2 < s)\, ds\, d\mathcal{B}_1(x_1)\, d\mathcal{B}_2(x_2),$$

where $\mathcal{B}_1$ and $\mathcal{B}_2$ are the two Wiener random measures; see [32]. We also need to define $\mathcal{Z}_{n2} := n^{1-H-h}\sum_{t=1}^n X_t u_t/\gamma$, $Z_{n2} := n^{-h}\sum_{i=1}^n X_i/\gamma$, and $c_1 := E(X^2\sigma(X))$. We are now ready to state and prove the following:

LEMMA 2.2. *Suppose* (1.1)–(1.3) *and Assumption* 1 *hold. Then,*

$$(2.5) \quad n^{1-H-h}\sum_{t=1}^{n} X_t\sigma_t u_t = \gamma c_1 \mathcal{Z}_{n2} + o_p(1),$$

$$(2.6) \quad n^{2-H-h}\bar{X}\bar{e} = \gamma\sigma_0 Z_{n1} Z_{n2} + o_p(1) \qquad \forall 1/2 < h, H < 1.$$

*Moreover, for* $H + h > 3/2$, $\mathrm{Correl}(\mathcal{Z}_{n2}, Z_{n1}Z_{n2})$ *converges to*

$$(2.7) \quad \frac{\sqrt{2(2H+2h-3)(2H+2h-2)}}{(2H+2h-1)}\sqrt{\frac{Hh}{(2H-1)(2h-1)}}.$$

PROOF. Let $H_j$ denote the $j$th Hermite polynomial, $j \geq 1$; see, for example, [32]. The Hermite expansion of the function $x\sigma(x)$ is $\sum_{j=0}^{\infty}\frac{c_j}{j!}H_j(x)$, where $c_j := E(X\sigma(X)H_j(X)), j \geq 0$. Since under Assumption 1, $c_0 = 0$, $c_1 = EX^2\sigma(X) \neq 0$, the Hermite rank $\min\{j \geq 1; c_j \neq 0\}$ of $x\sigma(x)$ is 1. Hence, in $L_2$,

$$(2.8) \quad n^{-1}\sum_{t=1}^{n} X_t\sigma_t u_t = \frac{c_1}{n}\sum_{t=1}^{n} X_t u_t + n^{-1}\sum_{t=1}^{n} u_t \sum_{j=2}^{\infty}\frac{c_j}{j!}H_j(X_t) =: S_n + T_n.$$

Using the independence of $X_i$'s and $u_i$'s, we obtain $\mathrm{Var}\, S_n = O(n^{-4+2H+2h})$. Because of the orthogonality of the Hermite polynomials,

$$\mathrm{Var}(T_n) = n^{-2}\sum_{s=1}^{n}\sum_{t=1}^{n} Eu_s u_t E\left\{\sum_{j=2}^{\infty}\frac{c_j}{j!}H_j(X_s)\sum_{k=2}^{\infty}\frac{c_k}{k!}H_k(X_t)\right\}$$

$$\leq Cn^{-2}\sum_{s=1}^{n}\sum_{t=1}^{n}|s-t|^{2H-2}|s-t|^{4h-4}$$

$$\leq Cn^{-6+2H+4h}\ln n = o(\mathrm{Var}\, S_n).$$

This fact and (2.8) readily yield (2.5). The claim (2.6) is proved similarly, using the fact that under Assumption 1, the Hermite rank of $\sigma(x)$ is 2.

To prove (2.7), let $\kappa_1 = G_X G_u \theta(H)\theta(h)$. By (2.1) and (2.8),

$$\mathrm{Var}\left(\sum_{t=1}^{n} X_t e_t/n\right) \sim c_1^2 \mathrm{Var}\left(\sum_{t=1}^{n} X_t u_t/n\right) \sim \frac{2c_1^2\kappa_1 n^{2H+2h-4}}{(2H+2h-3)(2H+2h-2)},$$

$$\mathrm{Var}(\bar{X}\bar{e}) \sim \sigma_0^2 \mathrm{Var}(\bar{X}\bar{u}) \sim \frac{\sigma_0^2\kappa_1 n^{2H+2h-4}}{(2H-1)(2h-1)Hh}.$$

Next, by the Hermite expansion of $x\sigma(x)$ and $\sigma(x)$, we obtain

$$E\left(n^{-1}\sum_{t=1}^{n} X_t e_t \bar{X}\bar{e}\right) \sim c_1\sigma_0 n^{-3}\sum_{t=1}^{n}\sum_{s=1}^{n}\sum_{k=1}^{n} E(X_t X_s)E(u_t u_k)$$

$$\sim c_1 \sigma_0 \kappa_1 n^{-3+2H+2h-4} \sum_{t=1}^{n} \sum_{s=1}^{n} \sum_{k=1}^{n} \left| \frac{t}{n} - \frac{s}{n} \right|^{2H-2} \left| \frac{t}{n} - \frac{k}{n} \right|^{2h-2}$$

$$\sim \frac{4\kappa_1 c_1 \sigma_0 n^{2H+2h-4}}{(2H-1)(2h-1)(2H+2h-1)}.$$

This proves (2.7). □

The following theorem gives a nonstandard second-order limiting distribution of $\hat{\beta}_1$ when $H \wedge h > 3/4$, where $\mathcal{Z}_2$ is as in (2.4) with $\mathcal{B}_1$ and $\mathcal{B}_2$ independent.

THEOREM 2.1. *Suppose* (1.1)–(1.3) *and Assumption 1 hold, with $\varepsilon_j$'s being $N(0,1)$ r.v.'s. Then, for $h \wedge H > 3/4$,*

(2.9) $$n^{2-H-h}(\hat{\beta}_1 - \beta_1) \to_d \gamma^{-1}[c_1 \mathcal{Z}_2 - \sigma_0 Z_1 Z_2].$$

PROOF. Lemma 2.2 combined with (2.3) yields that

$$n^{2-H-h}(\hat{\beta}_1 - \beta_1) = \frac{1}{\gamma}[c_1 \mathcal{Z}_{n2} - \sigma_0 Z_{n1} Z_{n2}] + o_p(1) \qquad \forall 1/2 < h, H < 1.$$

Using the derivations similar to those in the proof of Theorems 6.1 and 6.2 in [10], one verifies that $(\mathcal{Z}_{n2}, Z_{n1}, Z_{n2}) \to_d (\mathcal{Z}_2, Z_1, Z_2)$. Upon identifying $D_1, D_2$ there with $2(1-h), 2(1-H)$, respectively, one sees the condition $0 < D_1, D_2 < 1/2$ is equivalent to $h \wedge H > 3/4$. These facts together with (2.3) complete the proof of (2.9). □

Consistent estimates of $c_1, \sigma_0$ and $\gamma$ are $\sum_{i=1}^{n} X_i^2 V_i(X_i)/n$, $\sum_{i=1}^{n} V_i(X_i)/n$ and $s$, respectively, where $V_i$'s are defined in (5.2) below. However, the distribution of the limiting r.v. in (2.9) is not easy to determine and, hence, any decent inference about $\beta_1$ based on $\hat{\beta}_1$ appears to be infeasible in this case.

We shall next discuss asymptotic distribution of $\hat{\beta}_1$ when $u_t$'s form the moving average (1.2) and $H + h < 3/2$. This in turn is facilitated by the following lemma where $U_n := n^{-1/2} \sum_{t=1}^{n} \nu(X_t) u_t$ and $\gamma_\nu(k) = E\nu_0 \nu_k$.

LEMMA 2.3. *Suppose* (1.1)–(1.3) *holds. In addition, suppose $\nu$ is a measurable function such that $E\nu(X) = 0$, $E\nu^2(X) < \infty$, the Hermite rank of $\nu(X)$ is 1, and*

(2.10) $$\max_{\{0 \leq x \leq \ln n\}} |\nu(x)|/n^{1/2-\eta} \to 0, \qquad \text{for some } 0 < \eta < 1/2.$$

*Then, for $H + h < 3/2$, $U_n \to_d N(0, \kappa_2)$, where $\kappa_2^2 = \lim_{n \to \infty} EU_n^2 = \gamma_\nu(0)\gamma_u(0) + 2\lim_{n \to \infty} \sum_{k=1}^{n-1} \gamma_\nu(k)\gamma_u(k)$.*

PROOF. The proof uses the truncation method similar to the one used in [26]. The main idea here is to approximate $U_n$ by a weighted partial sum of the i.i.d. r.v.'s $\{\varepsilon_i\}$. Fix $H, h$ such that $H + h < 3/2$. Let $\nu_t := \nu(X_t)$ and $M = M_n > n^{(2h-1)/(2-2H)}$, and define

$$U_{n,M} := n^{-1/2} \sum_{t=1}^{n} \nu_t \sum_{j=-M}^{n} b_{t-j} \varepsilon_j.$$

Because the Hermite rank of $\nu_t$ is 1, by (A.20) below,

(2.11)
$$E(U_n - U_{n,M})^2 = n^{-1} \sum_{j=-\infty}^{-M-1} E\left(\sum_{t=1}^{n} \nu_t b_{t-j}\right)^2$$
$$\leq Cn^{2h-1} M^{-2+2H} \to 0.$$

Hence, it suffices to show that (a) $EU_{n,M}^2 \to \kappa_2^2$ and (b) $U_{n,M} \to_d N(0, \kappa_2)$.

Let $\mathcal{F} := \sigma\text{-}\{X_t, t \geq 1\}$. The claim (a) is implied by $E(U_{n,M}^2|\mathcal{F}) \to_p \kappa_2^2$. Let $d_{n,j} := \frac{1}{\sqrt{n}} \sum_{t=1}^{n} \nu_t b_{t-j}$. Then, $U_{n,M} = \sum_{j=-M}^{n} d_{n,j} \varepsilon_j$ and $E(U_{n,M}^2|\mathcal{F}) = \sum_{j=-M}^{n} d_{n,j}^2$. Recall $\gamma_u(k) = Eu_0 u_k = \sum_{j=0}^{\infty} b_j b_{j+k}$. Rewrite $E(U_{n,M}^2|\mathcal{F}) = A + 2B$, where

$$A = n^{-1} \sum_{j=-M}^{n} \sum_{t=1}^{n} \nu_t^2 b_{t-j}^2, \qquad B = n^{-1} \sum_{j=-M}^{n} \sum_{s=1}^{n-1} \sum_{t=s+1}^{n} \nu_t \nu_s b_{t-j} b_{s-j}.$$

But $A = A_1 + A_2$, where

$$A_1 = n^{-1} \sum_{j=-M}^{n} \sum_{t=1}^{n} (\nu_t^2 - \gamma_\nu(0)) b_{t-j}^2,$$
$$= n^{-1} \sum_{t=1}^{n} (\nu_t^2 - \gamma_\nu(0)) \left(\sum_{k=0}^{\infty} b_k^2 - \sum_{k=t+M+1}^{\infty} b_k^2\right) \xrightarrow{p} 0,$$
$$A_2 := \gamma_\nu(0) n^{-1} \sum_{j=-M}^{n} \sum_{t=1}^{n} b_{t-j}^2$$
$$= \gamma_\nu(0) n^{-1} \sum_{t=1}^{n} \left(\sum_{k=0}^{\infty} b_k^2 - \sum_{k=t+M+1}^{\infty} b_k^2\right) \to \gamma_\nu(0) \gamma_u(0),$$

because $\sum_{k=M}^{\infty} b_k^2 \to 0$, $\sum_{t=1}^{n} (\nu_t^2 - \gamma_\nu(0))/n \to 0$ and $\sum_{t=1}^{n} |\nu_t^2 - \gamma_\nu(0)|/n \to C < \infty$, a.s., by the Ergodic Theorem. Also, $\sup_n E|A_1| \leq C\gamma_\nu(0) < \infty$. Hence, $E|A_1| \to 0$ and $E|A - \gamma_\nu(0)\gamma_u(0)| \to 0$.

Next, let $g_j := n^{-1} \sum_{t=1}^{n-j} \nu_t \nu_{t+j}$. Then, one can rewrite $B = B_1 - B_2$, where $B_1 := \sum_{k=1}^{n-1} g_k \gamma_u(k)$ and $B_2 := \sum_{k=1}^{n-1} g_k \sum_{j=k+M+1}^{\infty} b_{k+j} b_j$. By (2.11),

$EB_2^2 \to 0$. For the term $B_1$, we have

$$B_1 = \sum_{k=1}^{n-1}[g_k - \gamma_\nu(k)]\gamma_u(k) + \sum_{k=1}^{n-1}\gamma_\nu(k)\gamma_u(k) =: B_{11} + B_{12} \quad \text{say.}$$

By applying Theorem 6 of [3] and the fact that the Hermite rank of the bivariate function $\nu_t\nu_{t+k} - \gamma_\nu(k)$ is 2, we obtain $\sup_k E|g_k - \gamma_\nu(k)| \leq Cn^{2h-2}$, for $1/2 < h < 1$, and hence, for $H + h < 3/2$,

$$E|B_{11}| \leq Cn^{2h-2}\sum_{k=1}^{n-1}|\gamma_u(k)| = O(n^{2H+2h-3}) \to 0.$$

Also, note that $\lim_n B_{12}$ exists for $H + h < 3/2$. These facts and the fact that $\kappa_2^2 = \gamma_\nu(0)\gamma_u(0) + 2\lim_n B_{12}$ complete the proof of (a).

The claim (b) is proved by showing that the conditional distribution of $U_{n,M}$, given $\mathcal{F}$, converges weakly to $N(0, \kappa_2)$. In view of the fact (a), by the Lindeberg–Feller theorem, this is equivalent to showing

$$(2.12) \qquad P\left(\max_{-M \leq j \leq n}|d_{n,j}| > \delta\right) \to 0 \qquad \text{for all } \delta > 0.$$

To prove this, recall from [5] that $\max_{1 \leq t \leq n}|X_t| = O_p(\ln n)$. By the Cauchy–Schwarz (C–S) inequality, we obtain that, for any integer $l > 0$,

$$\max_{-M \leq j \leq n}|d_{n,j}| = \max_{-M \leq j \leq n} n^{-1/2}\sum_{t=1}^{n}\nu_t b_{t-j}(I(|t-j| > l) + I(|t-j| \leq l))$$

$$\leq n^{-1/2}\left(\sum_{t=1}^{n}\nu_t^2\right)^{1/2}\max_{-M \leq j \leq n}\left(\sum_{t=1}^{n}b_{t-j}^2 I(|t-j| > l)\right)^{1/2}$$

$$+ n^{-1/2}\max_{1 \leq t \leq n}|\nu_t|\max_{-M \leq j \leq n}\sum_{t=1}^{n}|b_{t-j}|I(|t-j| \leq l)$$

$$= O_p\left(l^{-1+H} + n^{-1/2}\left(\max_{\{0 \leq x \leq \ln n\}}|\nu(x)|\right)l^{H-1/2}\right).$$

In view of (2.10), this upper bound is $o_p(1)$, for any $l = O(n^{2\eta/(2H-1)})$. Hence, (2.12) follows, thereby completing the proof of the lemma. $\square$

Now, take $\nu(x) = x\sigma(x)$ in the above lemma. Because $\sigma$ is an even function, the Hermite rank of $\nu(x)$ is 1. Also, the fact $\max_{1 \leq t \leq n}|X_t| = O_p(\ln n)$ and (2.6) imply $\bar{X}\bar{e} = o_p(n^{-1/2})$ for $H + h < 3/2$. Hence, we readily obtain the following:

THEOREM 2.2. *Suppose* (1.1) – (1.3), *and Assumption* 1 *hold. In addition, suppose* $EX^2\sigma^2(X) < \infty$ *and* (2.10) *holds with* $\nu(x) = x\sigma(x)$. *Then, for* $H + h < 3/2$, $n^{1/2}(\hat\beta_1 - \beta_1) \to_d N(0, \kappa_2/\gamma)$, where now $\kappa_2^2 = \lim_n \sum_{k=0}^n E\{X_0 \sigma(X_0) X_k \sigma(X_k)\} E(u_0 u_k)$.

An estimate of $\kappa_2$ is obtained as follows. Because, under Assumption 1, $|\text{Cov}(X_0\sigma(X_0)u_0, X_k\sigma(X_k)u_k)| \le Ck^{-4+2H+2h}$, for all sufficiently large $k$, the process $X_t\sigma(X_t)u_t$ is weakly dependent when $H + h < 3/2$. Thus, one may use the block bootstrap method to estimate $\kappa_2$ here using $X_t V_t(X_t)(Y_t - \hat\beta_0 - \hat\beta_1 X_t)$, $1 \le t \le n$, see [22], where $V_t$ is as in (5.2) below. Although we do not prove it here, such an estimator should be consistent for $\kappa_2$.

**3. Asymptotic distribution of $\hat\sigma^2(x)$.** In this section we shall investigate asymptotic distribution of the kernel type estimator $\hat\sigma^2(x)$ of $\sigma^2(x)$ in regression model (1.5) under the following:

ASSUMPTION 2. $E\|r(X)\|^2 < \infty$ and $A := Er(X)r(X)'$ is nonsingular.

An example of $r := (r_1, \ldots, r_q)$ satisfying this condition is $r_j(x) = x^j$, $j = 1, \ldots, q$.

To define $\hat\sigma^2(x)$, let $K$ be a density function on $[-1, 1]$, $b = b_n$ be sequence of positive numbers, $\phi$ denote the density of the $N(0,1)$ r.v., $\varphi(x) := \gamma^{-1}\phi((x-\mu)/\gamma)$, and $\varphi_n(x) := s^{-1}\phi((x-\bar X)/s)$. Let $K_b(x) \equiv K(x/b)/b$, $K_{bt}(x) := K_b(x - X_t)$, and define the kernel type estimator of $\sigma^2(x)$ to be

$$\hat\sigma^2(x) := \frac{1}{n\varphi_n(x)} \sum_{t=1}^n K_{bt}(x)\tilde e_t^2, \qquad \tilde e_t := Y_t - \tilde\beta' r(X_t).$$

Now fix an $x \in \mathbb{R}$ and consider the following additional assumptions.

ASSUMPTION 3. The density $K$ is symmetric around zero.

ASSUMPTION 4. The function $\sigma^2$ is twice continuously differentiable in a neighborhood of $x$.

ASSUMPTION 5. The bandwidth $b$ satisfies $b \to 0$, $n^{2h-1}(\ln n)^{-1} b \to \infty$ for $1/2 < h \le 3/4$, $n^{2-2h} b \to \infty$ for $3/4 < h < 1$.

To describe our results, we need to introduce $\mathcal{Z}_{n2}^* := n^{1-2H} \sum_{t=1}^n (u_t^2 - 1)$, $\mu_{r\sigma} := Er(X)\sigma(X)$. The proof of the following theorem is given in the Appendix. In it, $\mathcal{Z}_2^*$ is the $\mathcal{Z}_2$ of (2.4) with $\mathcal{B}_1 = \mathcal{B}_2$ and $\psi := (\psi_1^2 + \psi_2^2)^{1/2}$.

THEOREM 3.1. *Suppose* (1.2), (1.3), (1.5), *Assumptions* 2, 3 *and* 5 *hold and* $E\varepsilon^4 < \infty$.

(a) *In addition, suppose $x_1, \ldots, x_k$ are $k \geq 1$ points at which Assumption 4 holds and $r$ is continuous, $H < (1+h)/2$, and*

$$n^{1-h}b^2 \to 0. \tag{3.1}$$

*Then, $\{n^{1-h}(\hat{\sigma}^2(x_j) - \sigma^2(x_j)), j = 1, \ldots, k\}$ converges in distribution to $\{\frac{(x_j-\mu)}{\gamma}\sigma^2(x_j), j = 1, \ldots, k\}\psi Z$.*

(b) *In addition, suppose $x$ is a point at which $r$ is continuous and Assumption 4 holds, $H > (1+h)/2$ and*

$$n^{1-H}b \to 0. \tag{3.2}$$

*Then,*

$$\begin{aligned}
n^{2-2H}(\hat{\sigma}^2(x) - \sigma^2(x)) \\
= \sigma^2(x)\mathcal{Z}_{n2}^* + [\mu'_{r\sigma}A^{-1}r(x)r(x)'A^{-1}\mu_{r\sigma} \\
+ \mu'_{r\sigma}A^{-1}r(x)\sigma(x)]Z_{n1}^2 + o_p(1).
\end{aligned} \tag{3.3}$$

*Moreover, $(\mathcal{Z}_{n2}^*, Z_{n1}) \to_d (\mathcal{Z}_2^*, Z)$ and $\operatorname{Correl}(\mathcal{Z}_{n2}^*, Z_{n1}^2) \to \frac{2H}{4H-1}\sqrt{\frac{4H-3}{2H-1}}$.*

Consistent estimators of $\psi_1, \psi_2$ are obtained by plugging in the estimators of $H, G_u, G_X$ and $h$ in there while that of $\mu_{r\sigma}$ is $n^{-1}\sum_{i=1}^n r(X_i)V_i(X_i)$.

REMARK 3.1. Suppose we choose $b = O(n^{-\delta})$. Then Assumption 5 and (3.1) hold, for all $\delta$ in the range $(1-h)/2 < \delta < 2(1-h)$ whenever $h > 3/4$; and for all $\delta$ in the range $(1-h)/2 < \delta < 2h-1$, whenever $h \leq 3/4$ in the case (a). Similarly, in the case (b), Assumption 5 and (3.2) hold for $1-H < \delta < 2h-1$ whenever $h < 3/4$; and for $1-H < \delta < 2-2h$ whenever $h > 3/4$.

We also note here that by using the truncation method as in [2], the above Theorem 3.1 will continue to hold for a symmetric density kernel function $K$ with noncompact support and finite variance, for example, normal density.

**4. Estimation of $H$.** In this section we consider the problem of estimating $G_u, H$ in the model (1.5) based on $\tilde{e}_t := Y_t - \tilde{\beta}'r(X_t)$, $1 \leq t \leq n$.

For a process $\xi_t, 1 \leq t \leq n$, let $w_\xi(\lambda) := (2\pi n)^{-1/2}\sum_{t=1}^n \xi_t e^{it\lambda}$, $I_\xi(\lambda) := |w_\xi(\lambda)|^2$, $\lambda \in [-\pi, \pi]$, denote its discrete Fourier transform and periodogram, respectively, where $\mathbf{i} := (-1)^{1/2}$. Fix $1/2 < a_1 < a_2 < 1$. With $\lambda_j := 2\pi j/n$ and an integer $m \in [1, n/2)$, for $a_1 \leq \psi \leq a_2$, let

$$Q(\psi) := \frac{1}{m}\sum_{j=1}^m \lambda_j^{2\psi-1} I_{\tilde{e}}(\lambda_j), \qquad R(\psi) = \log Q(\psi) - (2\psi - 1)\sum_{j=1}^m \log \lambda_j.$$

Then the local Whittle estimators of $G_u$ and $H$ in the model (1.1) based on $\{\tilde{e}_t\}$ are defined to be $\hat{G}_u = Q(\hat{H})$, $\hat{H} = \arg\min_{\psi \in [a_1, a_2]} R(\psi)$, respectively.

The $\log(n)$ consistency of an analog of $\hat{H}$ and consistency of an analog of $\hat{G}_u$ in nonparametric homoscedastic regression models with $X_t = t/n$, $t = 1, \ldots, n$ is proved in [25]. The following theorem shows that these results continue to hold in the regression model (1.5) under much simpler restrictions on $m$ than those required in [25], partly due to the parametric nature of the model and partly due to random design.

THEOREM 4.1. *Suppose, in addition to* (1.2), (1.3), (1.5), *Assumptions* 2 *and* 5, *the following holds:*

$$(4.1) \qquad (\ln n)^4 \left( \left(\frac{m}{n}\right)^{2H-1} + \frac{m^{2(H-h)}}{n^{1+H-2h}} \right) \to 0.$$

*Then,* $\ln(n)(\hat{H} - H) \to_p 0$, $\hat{G} - G_u \to_p 0$.

The proof of this theorem is sketched in the Appendix. We note here that if $m = Cn^a$ for an $0 < a < 1$, then (4.1) holds for $H \geq h$. In the case $H < h$, it holds for any $a > (2h - H - 1)/(2h - 2H)$. In particular, in the case of Gaussian $\{u_t\}$'s, [14] shows that the optimal bandwidth $m$ equals $Cn^{4/5}$, which always satisfies (4.1).

**5. Regression model diagnostics.** In this section we investigate the weak convergence of $\widetilde{\mathcal{V}}_n$ under $H_0$ and the assumptions (1.2)–(1.3). The following Glivenko–Cantelli type result is used repeatedly in this connection: for a measurable real valued function $g$, with $E|g(X)| < \infty$,

$$(5.1) \qquad \sup_{x \in \mathbb{R}} \left| n^{-1} \sum_{t=1}^{n} g(X_t) I(X_t \leq x) - E g(X) I(X \leq x) \right| \xrightarrow{a.s} 0.$$

We are now ready to state and prove the following:

THEOREM 5.1. *Under* (1.2), (1.3), (1.5) *and Assumption* 2,

$$\sup_{x \in \bar{\mathbb{R}}} |n^{-H} \widetilde{\mathcal{V}}_n(x) - J_\sigma(x) Z_{n1}| = o_p(1).$$

*Hence, under* $H_0$, $n^{-H} \widetilde{\mathcal{V}}_n(x) \Longrightarrow J_\sigma(x) \psi_1 Z$, *in* $\mathcal{D}(\bar{\mathbb{R}})$, *and uniform metric.*

PROOF. Let $Z_n = \sum_{t=1}^{n} r(X_t) \sigma_t u_t$, $n\bar{A}_n := A_n$, and for an $x \in \bar{\mathbb{R}}$, let $\alpha(x) := Er(X)I(X \leq x)$, $L(x) := E\sigma^2(X)I(X \leq x)$,

$$\bar{\alpha}_n(x) := n^{-1} \sum_{t=1}^{n} r(X_t) I(X_t \leq x), \qquad F_\sigma(x) := E\sigma(X) I(X \leq x),$$

$$\mathcal{V}_n(x) := \sum_{t=1}^{n} \sigma_t u_t I(X_t \leq x), \qquad \mathcal{U}_n(x) := \sum_{t=1}^{n} u_t \{\sigma_t I(X_t \leq x) - F_\sigma(x)\}.$$

Now, assume $H_0$ holds. Using the Hermite expansion argument, we have $E(n^{-H}[\mathcal{U}_n(x) - \mathcal{U}_n(y)])^2 \leq Cn^{-2(1-h)}|L(y) - L(x)|$, for all $x,y \in \bar{\mathbb{R}}$. Then the chaining argument of [9] yields that $\sup_{x \in \mathbb{R}} |\mathcal{U}_n(x)| = o_p(1)$, and hence, $n^{-H}\mathcal{V}_n(x) = F_\sigma(x)Z_{n1} + u_p(1)$. By (2.1), we also have $n^{-H}Z_n = \mu_{r\sigma}Z_{n1} + o_p(1)$, $\bar{A}_n = A + o_p(1)$, and by (5.1), $\sup_{x \in \bar{\mathbb{R}}} \|\bar{\alpha}_n(x) - \alpha(x)\| = o_p(1)$. Note also that $J_\sigma(x) = F_\sigma(x) - \mu'_{r\sigma}A^{-1}\alpha(x)$. From these facts we readily obtain

$$n^{-H}\widetilde{\mathcal{V}}_n(x) := n^{-H}\sum_{t=1}^n [Y_t - \tilde{\beta}'r(X_t)]I(X_t \leq x)$$
$$= n^{-H}\mathcal{V}_n(x) - n^{-H}Z'_n\bar{A}_n^{-1}\bar{\alpha}_n(x) = J_\sigma(x)Z_{n1} + u_p(1).$$

This, uniform continuity of $J_\sigma$ and (2.2) complete the proof. $\square$

In order to implement the above result, we need a uniformly consistent estimator of $J_\sigma(x)$. One of the unknown entities in $J_\sigma$ is $\sigma(X)$. Because of unstable behavior of $\hat{\sigma}(X_t)$, we shall use an alternate estimator of $\sigma(x)$ based on the ideas of cross validation method that leave one observation out each time. For this purpose, we assume the design density is known, that is, $\mu, \gamma$ are known and take them to be $0, 1$ without the loss of generality. Let

$$(5.2) \quad \hat{\Lambda}_{-i}(x) := \left(\frac{1}{n-1}\sum_{t \neq i}^n K_{bt}(x)\tilde{e}_t^2\right)^{1/2}, \qquad V_i(x) := \hat{\Lambda}_{-i}(x)\phi^{-1/2}(x),$$

for $i = 1,\ldots,n$. Note that $V_i(x)$ is an estimator of $\sigma(x)$ and that of $J_\sigma(x)$ is

$$\hat{J}_n(x) = n^{-1}\sum_{t=1}^n V_t(X_t)I(X_t \leq x) - n^{-1}\sum_{t=1}^n r(X_t)V_t(X_t)\bar{A}_n^{-1}\bar{\alpha}_n(x).$$

To prove its uniform consistency, we need the following:

ASSUMPTION 6.  The function $\sigma$ has continuous first derivative.

THEOREM 5.2. *Suppose* (1.2), (1.3), (1.5), *and Assumptions* 2 *and* 6 *hold and that* $\mu = 0, \gamma = 1$. *In addition, suppose* $b \to 0$, $b^{-1}n^{2h-2} = O(1)$, $E\|r(2X)\|^2 < \infty$, $E\sigma^{2k}(X)r_j^4(X) < \infty$, *for* $j = 1,\ldots,q$, $k = 0,1$, *and* $E\sigma^2(X)\phi^{1/2}(X) < \infty$. *Then, under* $H_0$, $\sup_{x \in \mathbb{R}} |\hat{J}_n(x) - J_\sigma(x)| = o_p(1)$.

The proof of this theorem follows from the following lemma. Let $\tilde{\Lambda}_{-t}^2(x) = (n-1)^{-1}\sum_{i \neq t}^n K_{bi}(x)\sigma_i^2 u_i^2$.

LEMMA 5.1.  *Under the conditions of Theorem* 5.2,

$$(5.3) \qquad \max_{1 \leq t \leq n} E\{\tilde{\Lambda}_{-t}^2(X_t) - \sigma_t^2\phi(X_t)\}^2 \to 0,$$

$$(5.4) \qquad n^{-1}\sum_{t=1}^n |\hat{\Lambda}_{-t}(X_t) - \tilde{\Lambda}_{-t}(X_t)|^2\phi^{-1/2}(X_t) \xrightarrow{p} 0.$$

The proof of this lemma appears in the Appendix. We have the following:

COROLLARY 5.1. *Under the conditions of Theorem 5.2,*

$$\max_{1\le t\le n} E|\tilde{\Lambda}_{-t}(X_t) - \sigma_t \phi^{1/2}(X_t)|^4 \to 0, \tag{5.5}$$

$$n^{-1}\sum_{t=1}^{n} |[V_t(X_t) - \sigma_t]r_j(X_t)| \overset{p}{\to} 0, \qquad j = 1,\ldots,q. \tag{5.6}$$

PROOF. The claim (5.5) follows from (5.3) and the inequality $|a^{1/2} - b^{1/2}|^2 \le |a-b|, a \wedge b \ge 0$. We shall prove (5.6) for $j = 1$ only, it being similar for $j = 2, \ldots, q$. It suffices to show that

$$n^{-1}\sum_{t=1}^{n}|\hat{\Lambda}_{-t}(X_t) - \tilde{\Lambda}_{-t}(X_t)||r_1(X_t)|\phi^{-1/2}(X_t) = o_p(1), \tag{5.7}$$

$$n^{-1}\sum_{t=1}^{n}|\tilde{\Lambda}_{-t}(X_t)\phi^{-1/2}(X_t) - \sigma(X_t)||r_1(X_t)| = o_p(1). \tag{5.8}$$

By the Hölder inequality, the expectation of the l.h.s. of (5.8) is bounded above by $\sum_{t=1}^{n} E^{1/3}\{\tilde{\Lambda}_{-t}(X_t) - \sigma(X_t)\phi^{1/2}(X_t)\}^3 E^{2/3}\{r_1^{3/2}(X_t)\phi^{-3/4}(X_t)\}/n$. Since $E|r_1(X)|^{3/2}/\phi^{3/4}(X) = E|r_1(2X)|^{3/2} < \infty$, (5.8) follows from (5.5).

Next, by the C–S inequality, the l.h.s. of (5.7) is bounded above by

$$n^{-1}\left\{\sum_{t=1}^{n}|\hat{\Lambda}_{-t}(X_t) - \tilde{\Lambda}_{-t}(X_t)|^2\phi^{-1/2}(X_t) \cdot \sum_{t=1}^{n}r_1^2(X_t)\phi^{-1/2}(X_t)\right\}^{1/2}.$$

But because $n^{-1}\sum_{t=1}^{n} \frac{r_1^2(X_t)}{\sqrt{\phi(X_t)}} \to_{a.s.} E\frac{r_1^2(X)}{\sqrt{\phi(X)}} \le CEr_1^2(2X) < \infty$, (5.7) follows from this bound and (5.4). □

PROOF OF THEOREM 5.2. The proof follows from (5.6), the triangle inequality, the facts that $\|\bar{A}_n^{-1} - A^{-1}\| = o_p(1)$, and $\sup_x \|\bar{\alpha}_n(x) - \alpha(x)\| = o_p(1)$ implied by (5.1), in a routine fashion. □

A consequence of the above results is that whenever $\sup_x |J_\sigma(x)| \ne 0$, the test that rejects $H_0$, whenever, with $\hat{\psi}_1 := \hat{G}_u \theta(\hat{H})/\hat{H}(2\hat{H} - 1)$,

$$D_n := \frac{1}{n^{\hat{H}}\hat{\psi}_1 \sup_x |\hat{J}_n(x)|} \sup_x |\widetilde{\mathcal{V}}_n(x)| \ge z_{\alpha/2}, \tag{5.9}$$

is of the asymptotic size $\alpha$. Here $z_\alpha$ is the $100(1-\alpha)\%$ percentile of the $N(0,1)$ d.f. In the simple linear regression model with nonzero intercept, that is, when $r(x) = (1,x)'$, $J_\sigma(x) \equiv 0$ if and only if $\sigma(x)$ is constant in $x$.

In the case of a polynomial regression through the origin, $\sup_x |J_\sigma(x)| \neq 0$. In particular, the above test is applicable when fitting a heteroscedastic polynomial.

Proving consistency of the proposed test against a fixed alternative is a delicate matter. However, suppose $G_u, H$ and $\sigma$ are known such that $\sup_x |J_\sigma(x)| \neq 0$. Then the test that rejects $H_0$ whenever $\sup_x |\widetilde{\mathcal{V}}_n(x)| \geq n^H z_{\alpha/2} \psi_1 \sup_x |J_\sigma(x)|$ will be consistent against the alternative $\mu(x) = \beta' r(x) + \ell(x)$, for all $x$, where $\ell$ is such that $E\ell^2(X) < \infty$ and

$$\sup_{x \in \overline{\mathbb{R}}} |E[\ell(X) - E(r(X)\ell(X))'\bar{A}^{-1}r(X)]I(X \leq x)| \neq 0.$$

In the case these parameters are unknown, the above test (5.9) will be consistent against this alternative, provided estimators of these parameters continue to be consistent under the given alternative.

**6. Numerical results.** This section contains a simulation study and a real data application.

6.1. *A simulation study.* In this simulation we take $r(x) = (1, x)'$, $\beta_0 = 0$, $\beta_1 = 2$ and $\sigma^2(x) = 1 + x^2$. The errors $\{u_t\}$ are taken to be FARIMA$(0, H - 1/2, 0)$ with standardized Gaussian innovations and $\{X_t\}$ is taken to be fractional Gaussian noise with the LM parameter $h$. The values of $H, h$ range in the interval $[0.6, 0.95]$ with increments of 0.05. These processes were generated using the codes given in Chapter 12 of [4].

We first concentrate on the properties of $\hat{\beta}_1$ and $\hat{H}$. Table 1 provides the root mean square errors (RMSE) of the LSE $\hat{\beta}_1$ with sample size 500 and 2000 replications. As can be seen from this table, when $H + h$ increases, so does the RMSE of $\hat{\beta}_1$. Typically, when $H + h < 3/2$, the RMSE is small.

Table 2 provides the RMSE's of the local Whittle estimator $\hat{H}$ of $H$ based on $\hat{\varepsilon}_t = Y_t - \hat{\beta}_1 X_t$, $1 \leq t \leq 500$, repeated 1000 times. From this table, we observe that, for $H \leq 0.85$, the overall RMSE is less than 0.072 and stable regardless of the values of $h$.

TABLE 1
*RMSE of the LSE $\hat{\beta}_1$ for sample size $n = 500$*

| $H \setminus h$ | 0.60 | 0.65 | 0.70 | 0.75 | 0.80 | 0.85 | 0.90 | 0.95 |
|---|---|---|---|---|---|---|---|---|
| 0.60 | 0.00873 | 0.00864 | 0.00881 | 0.00980 | 0.01046 | 0.01151 | 0.01352 | 0.01922 |
| 0.65 | 0.00843 | 0.00950 | 0.01074 | 0.01170 | 0.01231 | 0.01345 | 0.01758 | 0.02473 |
| 0.70 | 0.01039 | 0.01011 | 0.01141 | 0.01354 | 0.01463 | 0.01760 | 0.02150 | 0.03410 |
| 0.75 | 0.01082 | 0.01212 | 0.01354 | 0.01545 | 0.01942 | 0.02273 | 0.03036 | 0.04653 |
| 0.80 | 0.01229 | 0.01395 | 0.01765 | 0.01924 | 0.02438 | 0.03328 | 0.04788 | 0.07352 |
| 0.85 | 0.01407 | 0.01776 | 0.02180 | 0.02828 | 0.03621 | 0.04875 | 0.07040 | 0.12540 |
| 0.90 | 0.01859 | 0.02369 | 0.03099 | 0.03984 | 0.05404 | 0.08341 | 0.12010 | 0.20873 |
| 0.95 | 0.02572 | 0.03406 | 0.05188 | 0.06472 | 0.11374 | 0.17622 | 0.27383 | 0.49624 |

TABLE 2
*RMSE of $\hat{H}$ based on $Y_t - \hat{\beta}_1 X_t$ for sample size $n = 500$*

| $H \setminus h$ | 0.60 | 0.65 | 0.70 | 0.75 | 0.80 | 0.85 | 0.90 | 0.95 |
|---|---|---|---|---|---|---|---|---|
| 0.60 | 0.03964 | 0.03867 | 0.03931 | 0.03870 | 0.03774 | 0.03940 | 0.03985 | 0.03938 |
| 0.65 | 0.04276 | 0.04273 | 0.04402 | 0.04261 | 0.04330 | 0.04511 | 0.04158 | 0.04317 |
| 0.70 | 0.04808 | 0.04750 | 0.04770 | 0.04819 | 0.05041 | 0.04856 | 0.04746 | 0.04690 |
| 0.75 | 0.05357 | 0.05478 | 0.05537 | 0.05473 | 0.05382 | 0.05075 | 0.05290 | 0.04940 |
| 0.80 | 0.06228 | 0.06303 | 0.06231 | 0.05918 | 0.05822 | 0.05973 | 0.05813 | 0.05509 |
| 0.85 | 0.07076 | 0.07202 | 0.07075 | 0.06699 | 0.06584 | 0.06310 | 0.06217 | 0.05773 |
| 0.90 | 0.08334 | 0.08323 | 0.08078 | 0.07803 | 0.076877 | 0.07237 | 0.06785 | 0.06508 |
| 0.95 | 0.11288 | 0.11096 | 0.10891 | 0.10300 | 0.09456 | 0.08519 | 0.07816 | 0.06718 |

TABLE 3
*Ranges for $\delta$ of the bandwidths for estimation $\sigma$*

| $H \setminus h$ | 0.65 | 0.75 | 0.85 | 0.95 |
|---|---|---|---|---|
| 0.65 | (a) (0.175, 0.3) | (a) (0.125, 0.5) | (a) (0.075, 0.3) | (a) (0.025, 0.1) |
| 0.75 | (a) (0.175, 0.3) | (a) (0.125, 0.5) | (a) (0.075, 0.3) | (a) (0.025, 0.1) |
| 0.85 | (b) (0.15, 0.3) | (a) (0.125, 0.5) | (a) (0.075, 0.3) | (a) (0.025, 0.1) |
| 0.95 | (b) (0.05, 0.3) | (b) (0.05, 0.5) | (b) (0.05, 0.3) | (a) (0.025, 0.1) |

TABLE 4
*Summary of $ASE(\hat{\sigma}^2)$ for $H = 0.65$*

| $h \setminus$ Summary | Bandwidth | Q1 | Median | Mean | Q3 |
|---|---|---|---|---|---|
| 0.65 | $3n^{-0.2}$ | 0.0261 | 0.0369 | 0.0424 | 0.0512 |
| 0.75 | $3.5n^{-0.2}$ | 0.0256 | 0.0383 | 0.0432 | 0.0557 |
| 0.85 | $4n^{-0.2}$ | 0.0273 | 0.0417 | 0.0595 | 0.0617 |
| 0.95 | $1.5n^{-0.099}$ | 0.0366 | 0.0663 | 0.1138 | 0.1058 |

TABLE 5
*Summary of $ASE(\hat{\sigma}^2)$ for $H = 0.75$*

| $h \setminus$ Summary | Bandwidth | Q1 | Median | Mean | Q3 |
|---|---|---|---|---|---|
| 0.65 | $4n^{-0.2}$ | 0.0442 | 0.0711 | 0.0887 | 0.1127 |
| 0.75 | $4n^{-0.2}$ | 0.0465 | 0.0652 | 0.0888 | 0.1076 |
| 0.85 | $4n^{-0.2}$ | 0.04667 | 0.0774 | 0.1043 | 0.1252 |
| 0.95 | $2n^{-0.099}$ | 0.0627 | 0.0995 | 0.2190 | 0.1902 |

Next, to assess the finite sample behavior of $\hat{\sigma}^2$, we simulated the estimator $\hat{\sigma}^2(x)$ for the values of $x$ in the grid $x_1 = -1.50, x_2 = -1.49, \ldots,$ $x_{301} = 1.50$, and for $0.65 \leq H, h \leq 0.95$. We used the built-in smoothing function of the R program with normal kernel and sample size 500 repeated

TABLE 6
*Summary of $ASE(\hat{\sigma}^2)$ for $H = 0.85$*

| $h \backslash$ Summary | Bandwidth | Q1 | Median | Mean | Q3 |
|---|---|---|---|---|---|
| 0.65 | $4.5n^{-0.2}$ | 0.1562 | 0.2724 | 0.5402 | 0.5584 |
| 0.75 | $6n^{-0.2}$ | 0.1594 | 0.3113 | 0.5449 | 0.6330 |
| 0.85 | $5n^{-0.2}$ | 0.1625 | 0.3252 | 0.5475 | 0.6103 |
| 0.95 | $2.5n^{-0.099}$ | 0.1704 | 0.3155 | 0.7092 | 0.6235 |

500 times. The ranges for $\delta$ in the bandwidths $b = Cn^{-\delta}$ are given in Table 3 according to the Remark 3.1. The symbols (a) and (b) indicate "Case a" and "Case b" of Theorem 3.1, respectively. Based on Table 3, for convenience, we used $\delta = 0.2, b = Cn^{-.2}$ in our simulations for all cases of $H$ and $h$ considered except when $h = 0.95$. In this case, we used $\delta = 0.099$. The constant $C$ is adjusted for different values of $H$ and $h$ according to the average squared errors: $ASE := \sum_{k=1}^{301}(\hat{\sigma}^2(x_k)/\sigma^2(x_k) - 1)^2/301$. We record those $C$ values which possibly make ASE the smallest. Some summary statistics of ASE are reported in Tables 4–7. It can be seen that the estimator $\hat{\sigma}^2(x)$ is relatively stable for the values of $H, h \leq 0.85$. Similar results are observed when we replace the normal kernel by the kernel function $K(x) = 0.5(1 + \cos(x\pi))I(|x| \leq 1)$ or the uniform kernel.

6.2. *Application to a foreign exchange data set.* In this section we shall apply the above proposed lack-of-fit test to fit a simple linear regression model with heteroscedastic errors to some currency exchange rate data obtained from www.federalreserve.gov/releases/H10/hist/. The data are noon buying rates in New York for cable transfers payable in foreign currencies. We use the currency exchange rates of the United Kingdom Pounds (UK£) vs. US$ and the Switzerland Franc (SZF) vs. US$ from January 4, 1971 to December 2, 2005. We first delete missing values and obtain 437 monthly observations. The symbols $X = $ dlUK and $Y = $ dlSZ stand for differenced log exchange rate of UK£vs. US$ and SZF vs. US$, respectively. We obtain

mean(dlUK) $= -0.0001775461$, Stdev(dlUK) $= 0.001701488$,

TABLE 7
*Summary of $ASE(\hat{\sigma}^2)$ for $H = 0.95$*

| $h \backslash$ Summary | Bandwidth | Q1 | Median | Mean | Q3 |
|---|---|---|---|---|---|
| 0.65 | $6n^{-0.2}$ | 1.153 | 3.214 | 16.24 | 11.83 |
| 0.75 | $7n^{-0.2}$ | 1.137 | 3.078 | 14.75 | 11.25 |
| 0.85 | $7.5n^{-0.2}$ | 1.018 | 2.611 | 12.77 | 11.59 |
| 0.95 | $4.5n^{-0.099}$ | 1.136 | 3.374 | 12.57 | 11.85 |

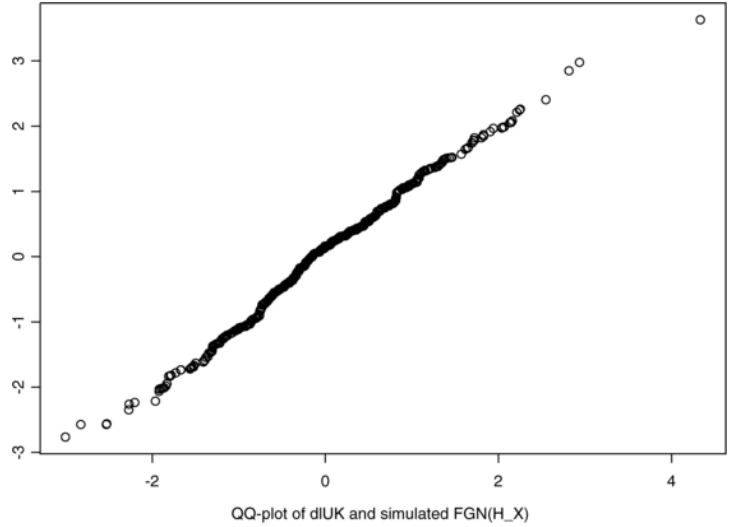

Fig. 1. *QQ-plot of dlUK.*

$$\text{mean(dlSZ)} = -0.00004525129, \qquad \text{Stdev(dlSZ)} = 0.001246904.$$

The local Whittle estimates of the LM parameters of dlUK and dlSZ processes, respectively, are 0.6610273 and 0.7147475. In computing these estimates we used $m = [n/8] = 54$.

Comparing the $X$-process with a simulated fractional Gaussian noise with $\hat{h} = 0.6610273$ and $n = 437$, Figure 1 suggests that the marginal distribution of $X$ is Gaussian.

Next, we regress $Y$ on $X$, using the normal density kernel regression function estimator and a simple linear regression model. Both of these estimates are depicted in Figure 2. They display a negative association between $X$ and $Y$. The estimated linear equation is $\hat{Y} = -0.000118775 - 0.4141107X$, with a residual standard error of 0.00102992.

Figure 3 provides the nonparametric kernel estimator of $\sigma(x)$ when regressing $Y$ on $X$ with $K(x) = 0.5(1 + \cos(x\pi))I(|x| \leq 1)$.

The estimators of $H$ based on $\hat{\varepsilon} = Y - \hat{\beta}X$ and $\hat{u} = (Y - \hat{\beta}X)/\hat{\sigma}(X)$ are equal to 0.6046235 and 0.6246576, respectively. This again suggests the presence of long memory in the error process.

Finally, to check if the regression of $Y$ on $X$ is simple linear, we obtain $D_n = 0.4137897$ with the asymptotic $p$-value 66%. As expected, this test fails to reject the null hypothesis that there exists a linear relationship between these two processes.

## APPENDIX

This section contains some preliminaries and proofs. To begin with we give a reduction principle involving the kernel function $K_b$. Let $\mathcal{G} := \{\nu : E\nu^2(X) < \infty\}$, $\mu = EX$, $\gamma^2 = \mathrm{Var}(X)$ and $Z = (X-\mu)/\gamma$. Then Hermite expansion of a $\nu \in \mathcal{G}$ is equal to $\sum_{j\geq 0}(c_j/j!)H_j(Z)$, where now $c_j = E\nu(\gamma Z + \mu)H_j(Z)$. We also need the fact that for any auto-covariance function $c(k) \sim Ck^{-2(1-\delta)}$, $k \to \infty$, $1/2 < \delta < 1$,

$$n^{-2}\sum_{t=1}^n \sum_{s\neq t} c^2(|t-s|) = O\left(\frac{1}{n}\right) + O\left(\frac{\log n}{n}\right) I(\delta = 3/4)$$

(A.1)
$$+ O(n^{4\delta - 4})I(3/4 < \delta < 1).$$

We are now ready to state and prove the following:

LEMMA A.1. *Let $X_t, t \in \mathbb{Z}$, be a stationary Gaussian process with $\mu = EX$ and $\gamma^2 = \mathrm{Var}(X)$. Let $\xi$ be a real valued measurable function on $\mathbb{R}$, $K$ be a density kernel on $\mathbb{R}$, and $b = b_n$ be a sequence of positive numbers, $b \to 0$. In addition, suppose $x \in \mathbb{R}$ is such that*

(A.2)
$$\sup_{n\geq 1} \int K^2(v)\xi^2(x-bv)\varphi(x-bv)\,dv < \infty.$$

*Let $\mu_b(x) = EK_b(x-X)\xi(X)$. Then,*

$$n^{-1}\sum_{t=1}^n (K_{bt}(x)\xi(X_t) - \mu_b(x)) - \xi(x)\frac{x-\mu}{\gamma}\varphi\left(\frac{x-\mu}{\gamma}\right)n^{-1}\sum_{t=1}^n \left(\frac{X_t-\mu}{\gamma}\right)$$

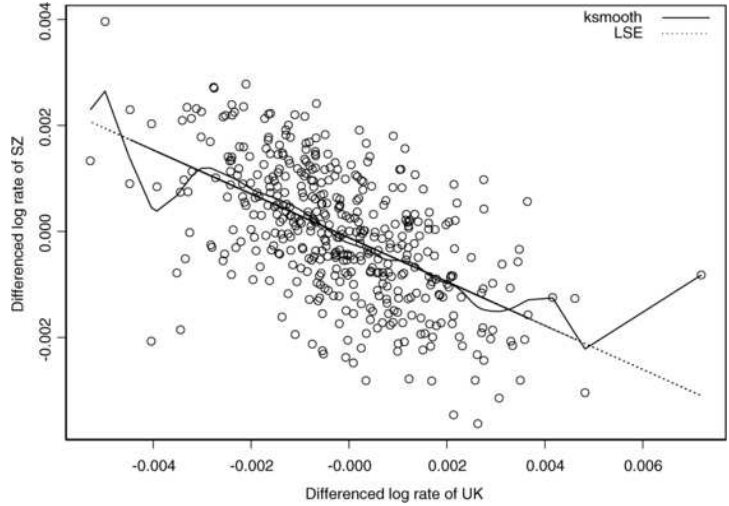

FIG. 2. *Kernel estimation of $r(x)$.*

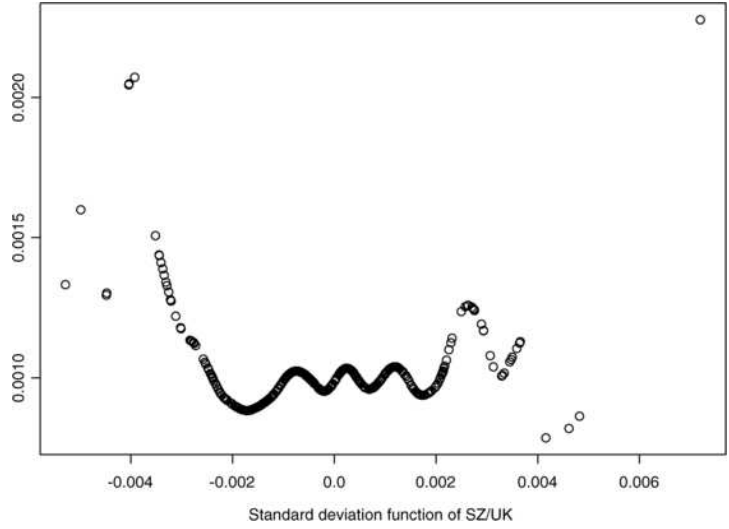

Fig. 3. *Kernel estimation of $\sigma(x)$.*

$$= O_p\Big(\frac{1}{\sqrt{nb}}\Big), \qquad\qquad 1/2 < h < 3/4,$$

$$= O_p\Big(\frac{1}{\sqrt{nb}}\Big) + O_p\Big(\Big(\frac{\log(n)}{nb}\Big)^{1/2}\Big), \qquad h = 3/4,$$

$$= O_p\Big(\frac{1}{\sqrt{nb}} + \frac{1}{\sqrt{bn^{4-4h}}}\Big), \qquad\qquad 3/4 < h < 1.$$

PROOF. Without loss of generality, assume $\mu = 0, \gamma = 1$. Let $x$ be as in (A.2). Let $\nu_n(X) := \sqrt{b}[K_b(x-X)\xi(X) - \mu_b(x)]$. For each $n \geq 1$, $E\nu_n(X) \equiv 0$, and $E\nu_n^2(X) \leq bEK_b^2(x-X)\xi^2(X) = \int K^2(v)\xi^2(x-bv)\varphi(x-bv)\,dv$. Hence, under (A.2), $\sup\{E\nu_n^2(X), n \geq 1\} < \infty$, so that $\nu_n(X) \in \mathcal{G}$, $\forall n \geq 1$. This in turn implies that $\sup_{n\geq 1}\sum_{j=1}^\infty c_{nj}^2/j! < \infty$, where, $\forall j \geq 1$,

$$c_{nj} := \sqrt{b}\int K_b(x-y)\xi(y)H_j(y)\varphi(y)\,dy = \sqrt{b}\{H_j(x)\xi(x)\varphi(x) + o(1)\}.$$

Hence,

$$E\Big\{\frac{1}{n}\sum_{t=1}^n[\nu_n(X_t) - c_{n1}X_t]\Big\}^2 = \frac{1}{n^2}\sum_{s=1}^n\sum_{t=1}^n\sum_{j=2}^\infty \frac{c_{nj}^2}{j!^2}EH_j(X_s)H_j(X_t)$$

$$\leq \sum_{j=2}^\infty \frac{c_{nj}^2}{j!}\Big\{\frac{1}{n} + \frac{1}{n^2}\sum_{t=1}^n\sum_{s\neq t}\gamma_X^2(|t-s|)\Big\}.$$

This and (A.1) applied to $c(k) = \gamma_X(k)$ complete the proof. □

We also need to recall the following result from [11]. Suppose $\{u_t\}$ is as in (1.2) with $E\varepsilon^4 < \infty$. Then, for all $1/2 < H < 1$,

(A.3) $\quad E(u_0^2 - 1)(u_t^2 - 1) = 2D^2 t^{2(2H-2)} + o(t^{2(2H-2)}), \quad t \to \infty,$

(A.4) $\quad Eu_0 u_s (u_t^2 - 1) \sim 2\gamma_u(t)\gamma_u(t-s), \quad |t-s| \to \infty,$

and for $s, t, r$ such that $|t-s|, |r-t|$ and $|s-r|$ all tending to infinity,

(A.5) $\quad Eu_0 u_s u_t u_r \sim \gamma_u(s)\gamma_u(r-t) + \gamma_u(t)\gamma_u(r-s) + \gamma_u(r)\gamma_u(t-s).$

Use (A.3)–(A.4), and argue as in the proof of Lemma A.1 with $\xi$ as in there, to obtain the following two facts under Assumptions 3 and 5, and (A.2):

(A.6) $\quad E\left\{n^{-1}\sum_{t=1}^{n}(K_{bt}(x)\xi(X_t) - \nu_b)u_t\right\}^2 = O\left(\frac{1}{nb}\right) + o\left(\frac{1}{n^{2-2H}}\right),$

(A.7) $\quad E\left\{n^{-1}\sum_{t=1}^{n}(K_{bt}(x)\xi(X_t) - \nu_b)(u_t^2 - 1)\right\}^2 = O\left(\frac{1}{nb}\right) + o\left(\frac{1}{n^{4-4H}}\right).$

Note that if $K$ is supported on $[-1,1]$ and $\xi$ is continuous at $x$, then (A.2) holds, for in this case the l.h.s. of (A.2) is bounded above by $C\max_{|z|\leq b}\xi^2(x-z) \to C\xi^2(x) < \infty$. It also holds if $\xi$ is bounded and $K$ is square integrable on $\mathbb{R}$. In particular, under Assumptions 3 and 4 and continuity of $r$ at $x$, it holds for $\xi(y) = \sigma^2(y)$ and $\xi(y) = \sigma(y)r(y)$, $y \in \mathbb{R}$.

Next, we give some inequalities that are useful in approximating an average of certain covariances of a square integrable function of a Gaussian vector by the corresponding average where the components of the Gaussian vector are i.i.d. Accordingly, let $E_k^0$ denote the expectation of a standard $k$-dimensional normal random vector. Let $A_{0,s,t}$ be the covariance matrix of $X_0, X_s, X_t$, $0 \leq s \leq t$, and $B_{0,s,t} = A_{0,s,t} - I_3 = ((b_{i,j}(s,t)))$, where $I_3$ is the $3 \times 3$ identity matrix. Let $\varrho_{s,t}$ denote the largest eigen value of $B_{0,s,t}$. From [23], Chapter 6.2, page 194, we obtain that $\varrho_{s,t} \leq \max_i \sum_{j=1}^{3}|b_{i,j}(s,t)|$. This in turn implies that

(A.8) $\quad \varrho_{s,t} \leq 3|\gamma_X(t-s)| \vee |\gamma_X(s)| \vee |\gamma_X(t)| \qquad \forall s \leq t, s, t \in \mathbb{Z}.$

For a square integrable function $g$ of $k$ r.v.'s, let $\|g\|_k^0 := (E_k^0 g^2)^{1/2}$ and $\tau_3(\tau_2)$ denote the Hermite rank of $g(X_0, X_1, X_2) - Eg(X_0, X_1, X_2)$ $(g(X_0, X_1) - Eg(X_0, X_1))$. Since both $\tau_3 \wedge \tau_2 \geq 1$, Theorem 2.1 of [27] yields that, for sufficiently large $|s-i|, |t-i|$ and $|t-s|$, $\exists$ a $C < \infty$ free of $i, s, t$, such that

(A.9) $\quad \begin{aligned}&|Eg(X_i, X_s, X_t) - E^0 g(X_i, X_s, X_t)| \\ &\leq C\|g\|_3^0 \varrho_{|s-i|,|t-i|}^{\tau_3/2} \leq C|\gamma_X(t-s)| \vee |\gamma_X(s-i)| \vee |\gamma_X(t-i)|,\end{aligned}$

$$|Eg(X_s, X_t) - E^0 g(X_s, X_t)|$$
(A.10)
$$\leq C\|g\|_2^0 \gamma_X^{\tau_2/2}(t-s).$$

In turn, (A.9) implies that, uniformly in $i = 1, \ldots, n$,

(A.11) $\quad (n-1)^{-2} \sum_{t \neq i} \sum_{s \neq i} (Eg(X_i, X_s, X_t) - E^0 g(X_i, X_s, X_t)) \to 0.$

PROOF OF THEOREM 3.1. Let $x$ be a point at which $r$ is continuous and Assumption 4 holds. Let $\tilde{d} := (\tilde{\beta} - \beta)$, $\tilde{\sigma}^2(x) := \sum_{t=1}^n K_{bt}(x) \sigma_t^2 u_t^2 / n\varphi(x)$,

$$S_n := \frac{1}{n\varphi(x)} \sum_{t=1}^n K_{bt}(x) \sigma_t r_t u_t, \qquad \Sigma_n := \frac{1}{n\varphi(x)} \sum_{t=1}^n K_{bt}(x) r_t r_t'.$$

Rewrite $\hat{\sigma}^2(x) - \sigma^2(x)$ as the sum

$$(\tilde{\sigma}^2(x) - \sigma^2(x)) + \frac{\varphi(x)}{\varphi_n(x)} [\tilde{d}' \Sigma_n \tilde{d} - 2\tilde{d}' S_n] + \left[\frac{\varphi(x)}{\varphi_n(x)} - 1\right] \tilde{\sigma}^2(x)$$

$$=: I^* + II^* + III^*.$$

Now let $\nu_b := EK_b(x - X)\sigma^2(X)$ and rewrite $I^* = I + II$, where

$$I := \frac{1}{n\varphi(x)} \sum_{t=1}^n K_{bt}(x) \sigma_t^2 (u_t^2 - 1), \qquad II = II_1 + II_2,$$

$$II_1 = \frac{1}{n\varphi(x)} \sum_{t=1}^n [K_{bt}(x) \sigma_t^2 - \nu_b], \qquad II_2 = \frac{\nu_b}{\varphi(x)} - \sigma^2(x).$$

First consider the term $II$. Use Assumption 4 to verify that

(A.12) $\qquad\qquad\qquad II_2 \leq Cb^2.$

Assumption 5 implies that $n^{2-2h} = o(nb)$. By Lemma A.1 applied to $\xi(y) = \sigma^2(y)$, we obtain that, under Assumptions 3–5,

(A.13) $\quad n^{1-h} II_1 = \gamma^{-1}(x - \mu) \sigma^2(x) Z_{n1} + o_p(1) \qquad \forall 1/2 < h < 1.$

Now consider the term $I$. Because $EI = 0$, (A.6)–(A.7) applied to $\xi(y) = \sigma^2(y)$ yield that $EI^2 = O((nb)^{-1})$, for $1/2 < H < 3/4$, $EI^2 = O((nb)^{-1} + n^{4H-4})$, for $3/4 < H < 1$, and $I = O_p((nb)^{-1/2} \ln^{1/2}(n))$, for $H = 3/4$. We summarize these results here: Under Assumptions 2–5, and $E\varepsilon^4 < \infty$,

$$I = O_p\left(\frac{1}{\sqrt{nb}} + n^{2H-2} I(H > 3/4) + \frac{\sqrt{\ln(n)}}{\sqrt{nb}} I(H = 3/4)\right),$$
(A.14)
$$II = O_p(n^{h-1} + b^2).$$

Next, consider $S_n$. By (A.3) with $\xi(y) = \sigma(y)r_j(y)$, $j = 1, \ldots, q$, we obtain

$$(A.15) \qquad n^{1-H}S_n = r(x)\sigma(x)n^{-H}\sum_{t=1}^{n} u_t + o_p(1).$$

Note also that $E(\Sigma_n) \to r(x)r'(x)$.

Lemma A.1 applied to $\xi(y) = (r(y)r'(y))_{i,j}$, $i, j = 1, \ldots, q$, yields

$$(A.16) \qquad \Sigma_n - E(\Sigma_n) = O_p(n^{h-1}).$$

To deal with the $III^*$ term, by (A.1) applied with $c(k) = \gamma_X(k)$, we have

$$s^2 - \gamma^2 = \frac{1}{n}\sum_{t=1}^{n}[(X_t - \mu)^2 - \gamma^2] - (\bar{X} - \mu)^2 = o_p(n^{h-1}).$$

This and the identity $s - \gamma = (s^2 - \gamma^2)/(s + \gamma)$ yield that $s - \gamma = o_p(n^{h-1})$. In turn, this fact, the continuity of $\varphi$ and the Taylor expansion yield that

$$(A.17) \qquad n^{1-h}\left(\frac{\varphi(x)}{\varphi_n(x)} - 1\right) = \frac{x - \mu}{\gamma}n^{1-h}(\bar{X} - \mu)/\gamma + o_p(1).$$

PROOF OF (a). Here $H < (1+h)/2$. From (A.12) and (A.14), one sees that in this case $n^{1-h}I = o_p(1) = n^{1-h}II_2$. Hence, by (A.13) and Assumption 5,

$$n^{1-h}I^* = \frac{x - \mu}{\gamma}\sigma^2(x)Z_{n1} + o_p(1).$$

By (A.15) and (A.16), $II^*$ is negligible, because

$$II^* = O_p(n^{2(H-1)} + n^{(2H-2)+(h-1)} + n^{2(H-1)}) = o_p(n^{h-1}).$$

These results, (A.17) and the fact that $\tilde{\sigma}^2(x) \to_p \sigma(x)$ yield

$$(A.18) \qquad n^{1-h}(\hat{\sigma}^2(x) - \sigma^2(x)) = \frac{x - \mu}{\gamma}\sigma^2(x)(Z_{n1} + Z_{n2}) + o_p(1),$$

where now $Z_{n2} = n^{-h}\sum_{t=1}^{n}(X_t - \mu)/\gamma$. By the independence of $X_t$'s and $u_t$'s, $Z_{n1}$ and $Z_{n2}$ are independent. Clearly, under the assumed conditions, (A.18) holds for each $x_1, \ldots, x_k$ given in part (a). Hence, part (a) follows from (2.2) and the Cramér–Wold device.

PROOF OF (b). In this case, $2H - 2 > h - 1$. Let $a = 2 - 2H$. Then, by (A.17), $n^a III^* = o_p(1)$. By (2.1) applied with $\nu(x) = r(x)\sigma(x)$,

$$(A.19) \qquad n^{1-H}\tilde{d} = n^{-H}A_n^{-1}\sum_{t=1}^{n} r_t\sigma_t u_t = A\mu_{r\sigma}Z_{n1} + o_p(1).$$

By (A.16), because $2H - 2 > h - 1$,

$$n^{2-2H}\tilde{d}'\Sigma_n\tilde{d} = n^{2-2H}\tilde{d}'[r(x)r(x)']\tilde{d} + o_p(1)$$
$$= \mu'_{r\sigma}A^{-1}[r(x)r(x)']A^{-1}\mu_{r\sigma}Z_{n1}^2 + o_p(1).$$

By (A.15) and (A.19), $n^{2-2H}\tilde{d}'S_n = \mu'_{r\sigma}A^{-1}r(x)\sigma(x)Z_{n1}^2 + o_p(1)$. Since $\varphi_n(x) \to_p \varphi(x) > 0$, we thus obtain

$$n^{2-2H}II^* = \mu'_{r\sigma}A^{-1}r(x)[r(x)'A^{-1}\mu_{r\sigma} + \sigma(x)]Z_{n1}^2 + o_p(1).$$

Next, consider $I^* = I + II$. Because $h > 1/2$, $H > (1+h)/2$ implies that $H > 3/4$. Hence, here the term $I$ is the dominating term. Recall $\mathcal{Z}_{n2}^* = n^{1-2H}\sum_{t=1}^n (u_t^2 - 1)$. By (A.7) applied with $\xi(y) = \sigma^2(y)$, we obtain

$$n^{2-2H}I^* = n^{2-2H}I + o_p(1) = \frac{\mu_b}{\varphi(x)}\mathcal{Z}_{n2}^* + o_p(1).$$

Upon combining these results, we obtain (3.3) in the case (b). The claim $(\mathcal{Z}_{n2}^*, Z_{n1}) \to_d (\mathcal{Z}_2^*, \psi_1 Z)$ is proved as in [28].

To prove the claim about $\text{Correl}(\mathcal{Z}_{n2}^*, Z_{n1}^2)$, note that $E(n^{1-H}\bar{u})^2 \to \psi_1^2$, and by (A.5),

$$E(\bar{u})^4 = \frac{4}{n^4}\sum_{t=1}^{n-1}\sum_{s_1=1}^{n-t}\sum_{s_2=1}^{n-t}\sum_{s_3=1}^{n-t} Eu_0 u_{s_1} u_{s_2} u_{s_3}$$

$$\sim \frac{4}{n^4}\sum_{t=1}^{n-1}\sum_{s_1=1}^{n-t}\sum_{s_2=1}^{n-t}\sum_{s_3=1}^{n-t}\{\gamma_u(s_1)\gamma_u(s_3-s_2)$$

$$+ \gamma_u(s_2)\gamma_u(s_1-s_3) + \gamma_u(s_3)\gamma_u(s_2-s_1)\}$$

$$\sim \frac{12G_X\theta(H)}{n^4}\sum_{t=1}^{n-1}\sum_{s_1=1}^{n-t}\gamma_u(s_1)\frac{(n-t)^{2H}}{H(2H-1)} \sim 3\psi_1^2 n^{4H-4}.$$

Hence, $\text{Var}(\bar{u}^2) \sim 2\psi_1^2 n^{-4+4H}$. By (A.3) and (A.4), similar calculations yield

$$E\left(\bar{u}^2 n^{-1}\sum_{t=1}^n (u_t^2 - 1)\right) \sim \frac{4G^2\theta^2(H)}{(2H-1)^2(4H-1)}n^{-4+4H},$$

$$\text{Var}\left(n^{-1}\sum_{t=1}^n (u_t^2 - 1)\right) \sim \frac{2G^2\theta^2(H)}{(4H-3)(2H-1)}n^{-4+4H}.$$

This completes the proof of Theorem 3.1. □

Next, to prove Theorem 4.1, we need the following preliminaries. Let $\xi$ be an arbitrary function such that $E\xi(X) = 0, E\xi^2(X) = 1$. Let $\xi_t := \xi(X_t)$ and $c_j := E\xi(\gamma Z + \mu)H_j(Z)$, $j \geq 1$. Let $\rho_k := \gamma_X(k)/\gamma^2$ and $\tau \geq 1$ be the Hermite rank of $\xi(X)$. Then, the auto-covariance function of the process $\xi$ is $\gamma_\xi(k) = \gamma_X(k)^\tau (\frac{c_\tau^2}{\tau!} + \sum_{j=\tau}^\infty \frac{c_{j+1}^2}{(j+1)!}\gamma_X(k)^{j+1-\tau})$, where the second term is bounded above by $\sum_{j\geq 1} c_j^2/j! = 1$. Therefore, there exists a constant $C = C(\tau, G_X)$ free of $k$, such that

(A.20) $$\gamma_\xi(k) \sim C\rho_k^\tau, \qquad k \to \infty.$$

LEMMA A.2. *With $\xi_t$ as above and $u_t$ as in* (1.2), *let $I_{\xi u}$ denote the periodogram of $\xi_t u_t$. Then, as $\lambda \to 0$,*

$$EI_{\xi u}(\lambda) = \begin{cases} O(\lambda^{\tau(2-2h)+1-2H}), & 0 < \tau(2-2h) + (2-2H) < 1; \\ O(\log n), & \tau(2-2h) + (2-2H) \geq 1. \end{cases}$$

PROOF. From (v.2.1), page 186 of [35] we obtain that, for any $0 < \alpha < 1$,

$$(A.21) \quad \sum_{t=1}^{n} t^{-\alpha} e^{\mathbf{i}t\lambda} \to \lambda^{\alpha-1} \Gamma(1-\alpha) \left( \sin \frac{\pi}{2}\alpha + \mathbf{i} \cos \frac{\pi}{2}\alpha \right), \qquad \lambda \to 0.$$

Since $\gamma_u(k) \sim Ck^{2H-2}$, by (1.4), (A.20), (A.21) with $0 < \alpha = \tau(2-2h) + (2-2H) < 1$ and the independence of $\xi_t$ and $u_t$, imply that, as $\lambda \to 0$,

$$EI_{\xi u}(\lambda) = \frac{1}{2\pi n} \sum_{j=1}^{n} \sum_{k=1}^{n} \gamma_\xi(j-k) \gamma_u(j-k) e^{\mathbf{i}(j-k)\lambda}$$

$$\sim Cn^{-1} \sum_{k=1}^{n} \sum_{t=k-n}^{k-1} t^{-[\tau(2-2h)+(2-2H)]} e^{\mathbf{i}t\lambda} \leq C\lambda^{\tau(2-2h)+1-2H}.$$

In the case $\tau(2-2h) + (2-2H) \geq 1$, $EI_{\xi u}(\lambda) \leq C \log n$, for all $\lambda \in [-\pi, \pi]$. □

The following fact proved in [7] is also needed for the proof of Theorem 4.1. Under (1.2) and when $m = o(n)$,

$$(A.22) \quad \sup_{0 < v \leq 1} [vm]^{-1} \sum_{j=1}^{[vm]} I_u(\lambda_j)/f_j \to 1 \qquad \text{a.s. } m \to \infty.$$

PROOF OF THEOREM 4.1. The basic proof is the same as in [25], with some difference in technical details. So we shall be brief, indicating only the main differences. With $\sigma_0 = E\sigma(X), r_0 := Er(X)$, let $\eta_t := e_t - \sigma_0 u_t$, $\xi := (\beta - \tilde{\beta})' r_0$ and $\zeta_t := (\beta - \tilde{\beta})' r(X_t)$. Then $\tilde{e}_t = \xi \zeta_t + \eta_t + \sigma_0 u_t$. Let $f_j = \lambda_j^{1-2H}$ and $\mathcal{D}_j =: [I_{\tilde{e}}(\lambda_j) - \sigma_0^2 I_u(\lambda_j)]/f_j$. According to the proof of Theorem 3 in [25], to prove Theorem 4.1 for $1/2 < H < 1$, it suffices to verify the following three claims:

$$(A.23) \quad \sum_{i=1}^{m-1} \left( \frac{i}{m} \right)^{2(a_1-H)+1} \frac{1}{i^2} \left| \sum_{j=1}^{i} \mathcal{D}_j \right| \xrightarrow{p} 0,$$

$$(A.24) \quad (\log n)^2 \sum_{i=1}^{m-1} \left( \frac{i}{m} \right)^{1-2\delta} \frac{1}{i^2} \left| \sum_{j=1}^{i} \mathcal{D}_j \right| \xrightarrow{p} 0 \qquad \text{for some small } \delta > 0,$$

(A.25)
$$\frac{(\log n)^2}{m} \sum_{j=1}^{m} \mathcal{D}_j \overset{p}{\to} 0.$$

We use the following elementary inequalities in verifying these conditions:

(A.26)
$$|I_{\tilde{e}}(\lambda) - \sigma_0^2 I_u(\lambda)| \leq 2\sigma_0 |I_u(\lambda) I_V(\lambda)|^{1/2} + I_V(\lambda),$$
$$I_V(\lambda) \leq 3(I_\xi(\lambda) + I_\zeta(\lambda) + I_\eta(\lambda)) \qquad \forall \lambda \in [-\pi, \pi],$$

where $V_t := \xi + \zeta_t + \eta_t$. Let $I_{Y,j} := I_Y(\lambda_j)$ for any process $Y_t$.

Recall that for the Dirichlet kernel $D_k(\lambda) := \sum_{t=1}^{k} e^{it\lambda}$, $|D_k(\lambda)| \leq C/\lambda$, for all $\lambda \in [-\pi, \pi], k \geq 1$. Also by (2.2) and (A.19),

(A.27)
$$n^{1-H} \|\beta - \tilde{\beta}\| = O_p(1).$$

These bounds imply that $\frac{I_{\xi,j}}{f_j} = O_p(n^{2H-3} \lambda_j^{2H-3})$, uniformly in $1 \leq j \leq m$.

Now, consider the terms $\zeta_t$. Assumption 2 implies that the Hermite ranks of $r_j(X) - r_{0j}, j = 1, 2, \ldots, q$, are at least one. Hence, by (A.26), (A.27) and (1.4), we obtain that, uniformly for $1 \leq j \leq m$,

(A.28)
$$\frac{I_{\zeta,j}}{f_j} = O_p(\lambda_j^{2(H-h)} n^{H-1}).$$

Similarly, using Lemma A.2 and the fact that the Hermite rank $\tau$ of $\sigma(X) - \sigma_0$ is at least 1, we obtain, uniformly for $1 \leq j \leq m$,

(A.29) $\frac{I_{\eta,j}}{f_j} = \begin{cases} O_p(\lambda_j^{\tau(2-2h)}), & 0 < \tau(2-2h) + (2-2H) < 1; \\ O_p(\lambda_j^{2H-1} \log n), & \tau(2-2h) + (2-2H) \geq 1. \end{cases}$

Now we are ready to verify (A.23)–(A.25). Let $\alpha_H := 2(a_1 - H)$. By changing the order of summation, the l.h.s. of (A.23) is bounded above by $Cm^{-\alpha_H - 1} \sum_{j=1}^{m} j^{\alpha_H} |\mathcal{D}_j|$, for $H > a_1$, and by $Cm^{-1} \log m \sum_{j=1}^{m} |\mathcal{D}_j|$, for $H = a_1$. But, by (A.26) and the C–S inequality,

$$\sum_{j=1}^{m} j^{\alpha_H} |\mathcal{D}_j| \leq C \sum_{j=1}^{m} j^{\alpha_H} \frac{\sigma_0 |I_{u,j} I_{V,j}|^{1/2} + I_{V,j}}{f_j}$$
$$\leq C \left[ \left( \sum_{j=1}^{m} j^{\alpha_H} \frac{|I_{u,j}|}{f_j} \sum_{j=1}^{m} j^{\alpha_H} \frac{|I_{V,j}|}{f_j} \right)^{1/2} + \sum_{j=1}^{m} j^{2\alpha_H} \frac{|I_{V,j}|}{f_j} \right].$$

We also have the following facts where $c_m := m^{-\alpha_H - 1}$:

$$\sum_{j=1}^{m} j^{a_H} \frac{I_{\xi,j}}{f_j} = \sum_{j=1}^{m} j^{c_m} \frac{I_{\eta,j}}{f_j} = \sum_{j=1}^{m} j^{c_m} \frac{I_{\zeta,j}}{f_j} = o_p(c_m^{-1}).$$

These bounds together with (A.22), (A.28) and (A.29) imply (A.23) for $H > a_1$. The proof of (A.23) for $H = a_1$ is similar. The conditions (A.24) and (A.25) are verified similarly. □

PROOF OF LEMMA 5.1. To prove (5.3), it suffices to show

$$\max_{1\le t\le n} E\left(\frac{1}{n-1}\sum_{j\ne t}^n K_{bj}(X_t)\sigma_j^2(u_j^2-1)\right)^2 \to 0, \tag{A.30}$$

$$\max_{1\le t\le n} E\left(\frac{1}{n-1}\sum_{j\ne t}^n K_{bj}(X_t)\sigma_j^2 - \sigma_t^2\phi(X_t)\right)^2 \to 0. \tag{A.31}$$

To prove (A.30), the expectation in the l.h.s. of (A.30) equals $A_{n,t} + B_{n,t}$, where $A_{n,t} := (n-1)^{-2}\sum_{j\ne t} E\{K_{bj}^2(X_t)\sigma_j^4\}E(u_j^2-1)^2$, and

$$B_{n,t} =: \frac{1}{(n-1)^2}\sum_{j\ne t}^n \sum_{k\ne t, k\ne j}^n E\{K_{bj}(X_t)K_{bk}(X_t)\sigma_j^2\sigma_k^2\}$$
$$\times E(u_j^2-1)(u_k^2-1).$$

Apply (A.9) and (A.11) to $g(X_t, X_j, X_k) = K_{bj}(X_t)K_{bk}(X_t)\sigma_j^2\sigma_k^2$. Note that for this $g$, $\|g\|_3^0 \le Cb^{-1}$. Hence, uniformly in $t$,

$$|B_{n,t}| \le \frac{C}{(n-1)^2}$$
$$\times \sum_{j\ne t}^n \sum_{k\ne t, k\ne j}^n \{E^0 K_{bj}(X_t)K_{bk}(X_t)\sigma_j^2\sigma_k^2$$
$$+ Cb^{-1}\varrho_{|t-j|,|t-k|}^{1/2}\}|E(u_j^2-1)(u_k^2-1)| \tag{A.32}$$
$$\le Cn^{4H-4} + Cb^{-1/2}n^{h-1+4H-4} \to 0,$$

by (A.3). Similarly, by (A.10), we obtain that, uniformly in $t$, $A_{n,t} \le C(nb)^{-1}$. Hence, (A.30) holds. To prove (A.31), rewrite the l.h.s. of (A.31) as

$$E(n-1)^{-2}\sum_{j,k\ne t}^n K_{bj}(X_t)K_{bk}(X_t)\sigma_j^2\sigma_k^2$$

$$- 2E(n-1)^{-1}\sigma^2(X_t)\phi(X_t)\sum_{j\ne t}^n K_{bj}(X_t)\sigma_j^2 + E\sigma^4(X)\phi^2(X)$$

$$=: C_{n,t} - 2D_{n,t} + E\sigma^4(X)\phi^2(X), \qquad \text{say}.$$

Similar to the argument in (A.32), by (A.9) and (A.10), the terms $C_{n,t}$ and $D_{n,t}$ tend to $E\sigma^4(X)\phi^2(X)$ uniformly in $t$, thereby proving (A.31). This completes the proof of (5.3).

Next consider (5.4). Arguing as for (5.5), it suffices to show that

$$n^{-1}\sum_{t=1}^n |\hat{\Lambda}_{-t}^2(X_t) - \tilde{\Lambda}_{-t}^2(X_t)|\phi^{-1/2}(X_t) \xrightarrow{p} 0. \tag{A.33}$$

But $|\hat{\Lambda}^2_{-t}(X_t) - \tilde{\Lambda}^2_{-t}(X_t)|$ is equal to

$$\left| \frac{1}{n-1} \sum_{i \neq t}^{n} K_{bi}(X_t)[(\beta - \hat{\beta})'r(X_i)]^2 + 2(\beta - \hat{\beta})'r(X_i)\sigma_i u_i \right|.$$

Moreover, (A.10) implies that, for $k = 0, 1, 2$, $1 \leq j \leq q$,

$$E \left| \frac{1}{n(n-1)} \sum_{t=1}^{n} \sum_{i \neq t}^{n} K_{bi}(X_t) \sigma_i u_i r_j^k(X_i) \phi^{-1/2}(X_t) \right|$$

$$\leq \frac{C}{n(n-1)} \sum_{t=1}^{n} \sum_{i \neq t}^{n} \int\int K_b(x-y)\sigma(y)|r_j(y)|^k \phi^{-1/2}(x) \phi_{i,t}(x,y) \, dx \, dy$$

$$\leq \frac{C}{n(n-1)} \sum_{t=1}^{n} \sum_{i \neq t}^{n} \int\int K_b(x-y)\sigma(y)|r_j(y)|^k \phi^{1/2}(x)\phi(y) \, dx \, dy$$

$$+ \frac{C}{n^{1-h}b^{1/2}} = O(1),$$

$$E \frac{1}{n(n-1)} \sum_{t=1}^{n} \sum_{i \neq t}^{n} K_{bi}(X_t)|r_j(X_i)|^k \phi^{-1/2}(X_t)$$

$$\leq \frac{C}{n(n-1)} \sum_{t=1}^{n} \sum_{i \neq t}^{n} \int\int K_b(x-y)|r_j(y)|^k \phi^{1/2}(x)\phi(y) \, dx \, dy + \frac{C}{n^{1-h}b^{1/2}}$$

$$= O(1).$$

In the above we used the assumptions $Er_j^4(X) < \infty$, $j = 1, \ldots, q$. Therefore, (A.33) holds because $\|\beta - \tilde{\beta}\| \to_p 0$. $\square$

**Acknowledgments.** The authors are grateful to the two referees for their constructive comments and to the Co-Editor Morris Eaton for his patience. They are especially grateful to the Associate Editor of this paper whose diligent and concerted effort helped to improve the presentation of the paper.

## REFERENCES

[1] AN, H. Z. and CHENG, B. (1991). A Kolmogorov–Smirnov type statistic with application to test for nonlinearity in time series. *Internat. Statist. Rev.* **59** 287–307.
[2] ANDREWS, D. W. K. (1995). Nonparametric kernel estimation for semiparametric models. *Econometric Theory* **11** 560–596. MR1349935
[3] ARCONES, M. (1994). Limit theorems for nonlinear functionals of a stationary Gaussian sequence of vectors. *Ann. Probab.* **22** 2242–2274. MR1331224
[4] BERAN, J. (1994). *Statistics for Long-Memory Processes*. Chapman and Hall, New York. MR1304490
[5] BERMAN, S. (1992). *Sojourns and Extremes of Stochastic Processes*. Wadsworth and Brooks/Cole, Pacific Grove, CA. MR1126464

[6] DAHLHAUS, R. (1995). Efficient location and regression estimation for long range dependent regression models. *Ann. Statist.* **23** 1029–1047. [MR1345213](MR1345213)

[7] DALLA, V., GIRAITIS, L. and HIDALGO, J. (2006). Consistent estimation of the memory parameter for nonlinear time series. *J. Time Ser. Anal.* **27** 211–251. [MR2235845](MR2235845)

[8] DAVYDOV, J. A. (1970). The invariance principle for stationary processes. *Theory Probab. Appl.* **15** 487–498. [MR0283872](MR0283872)

[9] DEHLING, H. and TAQQU, M. S. (1989). The empirical process of some long-range dependent sequences with an application to $U$-statistics. *Ann. Statist.* **17** 1767–1783. [MR1026312](MR1026312)

[10] FOX, R. and TAQQU, M. (1987). Multiple stochastic integrals with dependent integrators. *J. Multivariate Anal.* **21** 105–127. [MR0877845](MR0877845)

[11] GUO, H. and KOUL, H. L. (2007). Nonparametric regression with heteroscedastic long memory errors. *J. Statist. Plann. Inference* **137** 379–404. [MR2298945](MR2298945)

[12] GIRAITIS, L., KOUL, H. L. and SURGAILIS, D. (1996). Asymptotic normality of regression estimators with long memory errors. *Statist. Probab. Lett.* **29** 317–335. [MR1409327](MR1409327)

[13] HART, J. D. (1997). *Nonparametric Smoothing and Lack-of-Fit Tests*. Springer, New York. [MR1461272](MR1461272)

[14] HENRY, M. and ROBINSON, P. (1995). Bandwidth choice in Gaussian semiparametric estimation of long range dependence. *Athens Conference on Applied Probability and Time Series Analysis* **II** 220–232. *Lecture Notes in Statist.* **115**. Springer, New York. [MR1466748](MR1466748)

[15] HIDALGO, J. and ROBINSON, P. (2002). Adapting to unknown disturbance autocorrelation in regression with long memory. *Econometrica* **70** 1545–1581. [MR1929978](MR1929978)

[16] HO, H. and HSU, N. (2005). Polynomial trend regression with long memory errors. *J. Time Ser. Anal.* **26** 323–354. [MR2163286](MR2163286)

[17] KOUL, H. L. (1992). $M$-estimators in linear models with long range dependent errors. *Statist. Probab. Lett.* **14** 153–164. [MR1173413](MR1173413)

[18] KOUL, H. L. and MUKHERJEE, K. (1993). Asymptotics of $R$-, MD- and LAD-estimators in linear regression models with long range dependent errors. *Probab. Theory Related Fields* **95** 535–553. [MR1217450](MR1217450)

[19] KOUL, H. L. (1996). Asymptotics of $M$-estimators in non-linear regression with long-range dependent errors. *Athens Conference on Applied Probability and Time Series Analysis* **II**. *Lecture Notes in Statist.* **115** 272–290. Springer, New York. [MR1466752](MR1466752)

[20] KOUL, H. L., BAILLIE, R. and SURGAILIS, D. (2004). Regression model fitting with a long memory covariate process. *Econometric Theory* **20** 485–512. [MR2061725](MR2061725)

[21] KOUL, H. L. and STUTE, W. (1999). Nonparametric model checks for time series. *Ann. Statist.* **27** 204–237. [MR1701108](MR1701108)

[22] LAHIRI, S. N. (2003). *Reampling Methods for Dependent Data*. Springer, New York. [MR2001447](MR2001447)

[23] LUENBERGER, D. G. (1979). *Introduction to Dynamical Systems*: *Theory and Applications*. Wiley, New York.

[24] ROBINSON, P. M. (1995). Gaussian semiparametric estimation of long range dependence. *Ann. Statist.* **23** 1630–1661. [MR1370301](MR1370301)

[25] ROBINSON, P. M. (1997). Large-sample inference for nonparametric regression with dependent errors. *Ann. Statist.* **25** 2054–2083. [MR1474083](MR1474083)

[26] ROBINSON, P. and HIDALGO, J. (1997). Time series regression with long-range dependence. *Ann. Statist.* **25** 77–104. [MR1429918](MR1429918)

[27] SOULIER, P. (2001). Moment bounds and central limit theorem for functions of Gaussian vectors. *Statist. Probab. Lett.* **54** 193–203. MR1858634
[28] SURGAILIS, D. (1982). Zones of attraction of self-similar multiple integrals. *Lithuanian Math. J.* **22** 327–340. MR0684472
[29] STUTE, W. (1997). Nonparametric model checks for regression. *Ann. Statist.* **25** 613–641. MR1439316
[30] STUTE, W., THIES, S. and ZHU, L. X. (1998). Model checks for regression: An innovation process approach. *Ann. Statist.* **26** 1916–1934. MR1673284
[31] STUTE, W. and ZHU, L. X. (2002). Model checks for generalized linear models. *Scand. J. Statist.* **29** 333–354. MR1925573
[32] TAQQU, M. (1975). Weak convergence to fractional Brownian motion and to Rosenblatt process. *Z. Wahrsch. Verw. Gebiete* **31** 287–302. MR0400329
[33] YAJIMA, Y. (1988). On estimation of a regression model with long-memory stationary errors. *Ann. Statist.* **16** 791–807. MR0947579
[34] YAJIMA, Y. (1991). Asymptotic properties of the LSE in a regression model with long-memory stationary errors. *Ann. Statist.* **19** 158–177. MR1091844
[35] ZYGMUND, A. (2002). *Trigonometric Series*, 3rd ed. Cambridge Univ. Press. MR1963498

DEPARTMENT OF STATISTICS AND PROBABILITY
MICHIGAN STATE UNIVERSITY
EAST LANSING, MICHIGAN 48824-1027
USA
E-MAIL: guohongw@stt.msu.edu
koul@stt.msu.edu